\documentclass[12pt]{article}

\usepackage{amsmath}
\usepackage{wrapfig}
\usepackage{amssymb}
\usepackage{mathrsfs}
\usepackage{color}
\usepackage{graphicx}
\usepackage{mathabx}
\usepackage{dsfont}
\usepackage{float}
\usepackage{epsfig}
\usepackage{lscape}
\usepackage[bottom]{footmisc}
\usepackage{multirow}
\usepackage{booktabs}
\usepackage{caption}
\usepackage{natbib}

\newtheorem{theorem}{Theorem}[section]

\newtheorem{lemma}{Lemma}[section]
\newtheorem{remark}{Remark}[section]

\makeatletter

\@addtoreset{equation}{section}
\makeatother

\def\ind{ {{\rm 1}\hskip-2.2pt{\rm l}}}

\usepackage{natbib}

\begin{document}

\title{\bf Nonparametric estimation of conditional cure models for heavy-tailed distributions and under insufficient follow-up}

\author{
{\large Mikael E\textsc{scobar}-B\textsc{ach}}
\footnote{LAREMA, Universit\'e d'Angers, France. Email: \texttt{mikael.escobar-bach@univ-angers.fr}.  \vspace*{.1cm}} \\ {\it Universit\'e d'Angers} 
\and 
{\large Ingrid V\textsc{{an} K{eilegom}}} \footnote{ORSTAT, KU Leuven, Belgium. Email: \texttt{ingrid.vankeilegom@kuleuven.be}. Financial support from the European Research Council (2016-2021, Horizon 2020 / ERC grant agreement No.\ 694409) is  gratefully acknowledged.}\\ {\it KU Leuven}
}

\date{\today}

\maketitle

\begin{abstract}
When analyzing time-to-event data, it often happens that some subjects do not experience the event of interest.    Survival models that take this feature into account (called `cure models') have been developed in the presence of covariates. However, the current literature on nonparametric cure models with covariates cannot be applied when the follow-up is insufficient, i.e., when the right endpoint of the support of the censoring time is strictly smaller than that of the survival time of the susceptible subjects. In this paper we attempt to fill this gap in the literature by proposing new estimators of the conditional cure rate and the conditional survival function using extrapolation techniques coming from extreme value theory. We establish the asymptotic normality of the proposed estimators, and show how the estimators work for small samples by means of a simulation study.  We also illustrate their practical applicability through the analysis of data on the survival of colon cancer patients. 
\end{abstract}

\newpage

\section{Introduction} \label{section_introduction}

In survival analysis it often happens that the event of interest never occurs for a fraction of the subjects under study.  This is often the case in medical studies in which one is interested in the survival time of patients receiving a treatment for a certain disease.  If a patient gets cured from his/her disease, we will never observe his/her time to death for that particular disease.   Other common examples can be found in demography (time to marriage, time to first child), economics (time to finding a job), marketing (time to buying a product), education (time to learning how to do a certain task), etc.   A large body of papers have been published over the last 20 years in which survival models have been extended to take the cure fraction into account, and these models have (quite naturally) been called {\it cure models}.   The current literature includes a wide range of parametric, semiparametric and fully nonparametric models.  The literature on nonparametric cure models is rather scarce compared to the rich literature on the parametric and semiparametric counterpart. Among them, we can cite \cite{Xu2014}, \cite{Lopez2017}, \cite{Lopez2017b} and \cite{Chown2018}, who all consider covariates in the model. We refer to Maller and Zhou (1996) for a book-long overview of early references on cure models, and to \cite{Peng2014} and \cite{Amico2018} for recent review papers on cure models.  Cure models have been studied from different angles, including theoretical, methodological, modeling, computational and applied points of view.   The quantities of interest in these models are often the cure rate and the survival function of the susceptible or uncured individuals, i.e.\ those who will sooner or later experience the event of interest.  The latter two quantities are often allowed to depend on covariates, and many models (in particular the so-called mixture cure models) allow these two quantities to depend on different sets of covariates.   

When the survival time is subject to random right censoring, as is common in survival analysis, 
all cured subjects will be censored, whereas the non-cured ones can be either censored or uncensored.  Hence, it is clear that in order to identify the cure fraction we need to impose certain assumptions on the model.   A common way to identify a cure model is to impose the assumption of {\it sufficient follow-up}, which means that the right endpoint of the support of the censoring time is larger than the right endpoint of the support of the survival time of the non-cured subjects (conditional on the covariates in case there are covariates in the model).   See e.g.\ \citet{Maller1992} for more details.   In non- and semiparametric cure models, this assumption is standard.  An informal way to verify this assumption in practice is to check whether the \cite{Kaplan1958} estimator has a sufficiently long plateau, that contains several censored observations.  This is however a vague statement, and its judgement is rather subjective.  When the follow-up period is erroneously believed to be sufficiently long, the cure rate will be overestimated leading to possibly false (and too positive) conclusions regarding the effect of a treatment or a drug.   Therefore, it is important to have a method at hand that is able to correctly estimate the cure rate (and also the survival function of the uncured subjects), when the follow-up period is (possibly) insufficient.   

The goal of this paper is to provide a nonparametric method in the presence of covariates that allows to extrapolate the conditional survival function beyond the last data point in order to correctly estimate both the conditional cure rate and the survival function itself.  This will be achieved by using techniques from extreme value theory, assuming that the conditional distribution of the uncured subjects is heavy-tailed.   The use of extreme value theory in survival analysis is not new.  See e.g.\ \cite{Beirlant2001}, \cite{Einmahl2008}, \cite{Beirlant2010}, \cite{Gomes2011}, \cite{Worms2014} and \cite{Stupfler2016}, among others.  However, all these papers focus on the estimation of the survival function in the case where there is no cure fraction, which is intrinsically a much easier problem, as the identifiability of the cure fraction is not an issue in that case.   We will show that extreme value theory allows to construct an estimator of the conditional cure rate that converges to the correct cure rate when the sample size and the endpoint of the censoring distribution converge to infinity.

In the absence of covariates, \cite{Escobar2019} proposed a method to improve the estimation of the cure rate by using extreme value theory.  With respect to the latter paper, this paper makes several steps forward.  First, covariates are included in the model, which is an important improvement towards the application of the method in practice, but which complicates the theory considerably.  And second, in this paper we do not only estimate the conditional cure rate, but we also propose an estimator of the conditional extreme value index and of the conditional survival function.  The latter estimator will be based on the local Kaplan-Meier estimator proposed by \cite{Beran1981}, to which a term is added that corrects the bias caused by the insufficient follow-up.  Also note that the estimator that we propose for the conditional cure rate, does not reduce to the estimator proposed in \cite{Escobar2019} in the absence of covariates. By extending the degree of freedom, we propose an improved estimator that eases the estimation accuracy and the variability reduction in practice. This also comes with a cross-validation procedure that all together offers a fully data-driven solution.\\   

The remainder of the paper is organized as follows. In the next section we introduce our proposed estimators of the conditional extreme value index, the cure rate and the survival function.  In Section \ref{section_asymptotics} we first develop the asymptotic properties of the Beran estimator upto (and including) the right endpoint of the censoring distribution, and next we consider the asymptotics of the proposed estimators.  In particular, we estimate the survival function beyond this right endpoint. The finite sample performance of the proposed estimators is studied in Section \ref{section_simulation} for various models within a Fr\'echet domain of attraction, while in Section \ref{section_data} the method is applied on data coming from a study on the survival of colon cancer patients. Finally, the proofs of the main asymptotic results are relegated to the Appendix.

\section{The extrapolation method} \label{section_extrapolation}

\subsection{Notations and definitions}

We start by introducing some notations.  Let $Y$ denote the survival time of a subject, and let $X$ be a $p$-dimensional vector of covariates with domain $S_X \subset \mathbb{R}^p$.   We denote the conditional cure rate for a given vector $x \in S_X$ by 
\begin{eqnarray*}
1-p(x) =\mathbb{P}(Y=+\infty|X=x).
\end{eqnarray*}
In the presence of random right censoring, we do not always observe $Y$ but instead the observed variables are $T=\min(Y,C)$ and the censoring indicator $\delta=\ind_{\{Y \le C\}}$, where $C$ refers to a random censoring time that is assumed to be finite. In this context, all cured subjects, i.e.\ those subjects for which $Y$ is infinite, are censored, and among the non-cured or susceptible subjects, some or censored and others are not. The conditional sub-distribution $F$ of $Y$ can be written as
\begin{eqnarray} 
\label{model}
F(t|x) = \mathbb{P}(Y\leq t|X=x) = p(x) F_0(t|x),
\end{eqnarray}
where $F_0(\cdot|x)$ is the (proper) conditional distribution of the survival time for the susceptible subjects.  Model (\ref{model}) is called a {\it mixture cure model}, as the survival function can be written as $S(t|x) = 1-F(t|x) = 1-p(x) +p(x) S_0(t|x)$ (with $S_0=1-F_0$), which is a mixture of the survival function of the cured and non-cured sub-populations.  The conditional distribution of the censoring time is denoted by $G(t|x)=\mathbb{P}(C\leq t|X=x)$. The right endpoint of the support of the distribution $G(\cdot|x)$ is denoted by $\tau_c(x)$.

We will work under minimal conditions on the distributions, though we have to impose the usual identification assumption that $Y$ and $C$ are independent conditionally on $X$, which implies that $H(t|x)=\mathbb{P}(T \le t|X=x)$ satisfies $1-H(t|x)=(1-F(t|x))(1-G(t|x))$. In the sequel, we will also use the notations $H^u(t|x)=\mathbb{P}(T\leq t,\delta=1|X=x)=\int_{-\infty}^t (1-G(s^-|x)) dF(s|x)$ for the sub-distribution function of the uncensored observations and $\Lambda(\cdot|x)$ for the cumulative hazard function given by
\begin{eqnarray*}
\Lambda(t|x)=\int_{-\infty}^t \dfrac{dH^u(s|x)}{1-H(s^-|x)}.
\end{eqnarray*} 

\subsection{The main idea}

In the context of insufficient follow-up, part of the support of the survival time is not observable.  In the case of heavy tailed distributions, this simply reduces to $\tau_c(x)<+\infty$ since these distributions admit infinite right endpoints. One way to handle this problem is to assume that we dispose of some information on the queue of the survival function, that can be used to estimate the non-observable part.  This is the case in extreme value theory, where the tail information is characterized by a central parameter called the {\it extreme value index} (see e.g. \citealt{deHaan2006} for further details). We thus assume that $F_0(\cdot|x)$ is heavy-tailed, which means that it belongs to the Fr\'echet domain of attraction with conditional extreme value index $\gamma(x)>0$. This means that for any $y>0$,
\begin{eqnarray}
\label{doa}
\lim_{t\rightarrow+\infty}\dfrac{1-F_0 \big(t(1+y\gamma(x))|x\big)}{1-F_0(t|x)}=(1+y\gamma(x))^{-1/\gamma(x)}.
\end{eqnarray}
The extreme value index characterizes the family of all possible limiting distributions for large observations in a sample drawn from $F_0(\cdot|x)$ and drives its tail behaviour. In particular, it allows us to extrapolate the survival function beyond the support of the censoring distribution. This initial condition is similar to what has been done in \cite{Escobar2019}, although the latter paper does not consider covariates. To go further, we consider in that order the estimation of the conditional extreme value index, the conditional cure rate and the conditional survival function. 
To explain the idea, we replace $1+y\gamma(x)$ in (\ref{doa}) by $y$ and obtain 
\begin{eqnarray*}
\lim_{t\to+\infty}\dfrac{1-F_0(yt|x)}{1-F_0(t|x)}=y^{-1/\gamma(x)}.
\end{eqnarray*}
Since no observations lie outside of the censoring support, we replace $t$ by $\tau_c(x)$, that we assume large enough. By (\ref{model}), it turns out that
\begin{eqnarray*}
\dfrac{F_0(\tau_c(x)|x)-F_0(y\tau_c(x)|x)}{1-F_0(\tau_c(x)|x)}=\dfrac{F(\tau_c(x)|x)-F(y\tau_c(x)|x)}{p(x)-F(\tau_c(x)|x)}\simeq y^{-1/\gamma(x)}-1.\\
\end{eqnarray*}

Using at the same time $y^2$ and $y$ firstly allow us to obtain an expression that is free of the cure rate, namely
\begin{eqnarray*}
\dfrac{F(y^2\tau_c(x)|x)-F(y\tau_c(x)|x)}{F(y\tau_c(x)|x)-F(\tau_c(x)|x)}&=&\dfrac{F(y^2\tau_c(x)|x)-F(y\tau_c(x)|x)}{p(x)-F(\tau_c(x)|x)}\dfrac{p(x)-F(\tau_c(x)|x)}{F(y\tau_c(x)|x)-F(\tau_c(x)|x)}\\[.2cm]
&\simeq&\dfrac{y^{-2/\gamma(x)}-y^{-1/\gamma(x)}}{y^{-1/\gamma(x)}-1}=y^{-1/\gamma(x)},
\end{eqnarray*}
which gives us a first approximation for the conditional extreme value index, namely
\begin{eqnarray}
\label{evi}
\gamma(x)\simeq -\log(y)\Big/\log\left(\dfrac{F(y^2\tau_c(x)|x)-F(y\tau_c(x)|x)}{F(y\tau_c(x)|x)-F(\tau_c(x)|x)}\right).
\end{eqnarray}
The same argument gives us a approximation for the non-cure rate given by 
\begin{eqnarray}
\label{cr}
p(x) \simeq F(\tau_c(x)|x)+\dfrac{F(\tau_c(x)|x)-F(y\tau_c(x)|x)}{y^{-1/\gamma(x)}-1}.
\end{eqnarray}
Following the same ideas, we also propose to extrapolate $F_0(t|x)$ for $t>\tau_c(x)$. To do so, we obtain with $y=t/\tau_c(x)$,
\begin{eqnarray*}
\dfrac{F(\tau_c(x)|x)-F(t|x)}{p(x)-F(\tau_c(x)|x)}\simeq \left[\dfrac{t}{\tau_c(x)}\right]^{-1/\gamma(x)}-1,
\end{eqnarray*}
yielding 
\begin{eqnarray}
\label{F0}
F(t|x)\simeq F(\tau_c(x)|x) - \Big(p(x)-F(\tau_c(x)|x)\Big)\left(\left[\dfrac{t}{\tau_c(x)}\right]^{-1/\gamma(x)}-1\right).
\end{eqnarray}
In the next subsection, formulas (\ref{cr}), (\ref{evi}) and (\ref{F0}) will form the basis for the construction of our estimators of $\gamma(x)$, $p(x)$ and $F(t|x)$.

\subsection{The proposed estimators} 

Let $\{(T_i,\delta_i,X_i)\}_{1\leq i\leq n}$ be an i.i.d.\ sample drawn from the triplet $(T,\delta,X)$, denote by $T_{(i)}$ the $i$-th order statistic in the sample, and let $\delta_{(i)}$ and $X_{(i)}$ be the corresponding censoring indicator and covariate vector. The \cite{Beran1981} estimator is then given by
\begin{eqnarray*}
F_n(t|x)=1-\prod_{T_{(i)}\leq t}\left(1-\dfrac{W_h(x-X_{(i)})}{1-\sum_{j=1}^{i-1}W_h(x-X_{(j)})}\right)^{\delta_{(i)}},
\end{eqnarray*}
where for any $i=1,\ldots,n$,
\begin{eqnarray*}
W_h(x-X_i)=\dfrac{K_h(x-X_i)}{\sum_{j=1}^nK_h(x-X_j)},
\end{eqnarray*}
and $K_h(\cdot)=K(\cdot/h)/h^p$ with $K$ a kernel function and $h=h_n$ a non-random positive sequence such that $h_n\rightarrow 0$ as $n\rightarrow \infty$.\\

For the construction of the estimators of $\gamma(x)$, $p(x)$ and $F(t|x)$, we need to assume that the right endpoint $\tau_c(x)$ is the same for all $x$, i.e.\ $\tau_c(x)=\tau_c$.  A natural estimator of $\tau_c$ is given by $\tau_n=\max\{T_i : i=1,\ldots,n\}$.  The estimators of the conditional extreme value index $\gamma(x)$ and the non-cure rate $p(x)$ are then respectively given by
\begin{eqnarray} \label{gammahat}
\widehat\gamma(x)=-1/\log_{y_2}\left(\dfrac{F_n(y_2^2\tau_n|x)-F_n(y_2\tau_n|x)}{F_n(y_2\tau_n|x)-F_n(\tau_n|x)}\right)
\end{eqnarray}
and
\begin{eqnarray} \label{phat}
\widehat p(x)= F_n(\tau_n|x)+\dfrac{F_n(\tau_n|x)-F_n(y_1\tau_n|x)}{y_1^{-1/\widehat\gamma(x)}-1},
\end{eqnarray}
where $y_1,y_2\in(0,1)$ are tuning parameters. Note here that we have two different parameters against one as in \cite{Escobar2019}. This allows to compute $\widehat{\gamma}(x)$ separately with a dedicated parameter before estimating $\widehat{p}(x)$. Furthermore, this flexibility will be necessary to propose an efficient parameters selection in practice. 
Finally, for any $t\in\mathbb{R}$, the estimator of the distribution function $F(t|x)$ is defined by
\begin{eqnarray} \label{Fhat}
\widehat{F}(t|x)=F_n(t\wedge\tau_n|x)+\Big(\widehat p(x)-F_n(\tau_n|x)\Big)\left(1-\left[\frac{t}{\tau_n}\vee 1\right]^{-1/\widehat\gamma(x)}\right).
\end{eqnarray}
As will be shown in the next section, the estimators given in (\ref{gammahat}), (\ref{phat}) and (\ref{Fhat}) do not consistently estimate $\gamma(x)$, $p(x)$ and $F(t|x)$ for $t>\tau_c$, but instead they converge respectively to the functions (of $\tau_c$) given in (\ref{evi}), (\ref{cr}) and (\ref{F0}).  The latter quantities converge to the targeted quantities as $\tau_c$ tends to $\infty$. We denote these functions by
\begin{eqnarray*}
\gamma_{y_2,\tau_c}(x)&=&-1 \Big/ \log_{y_2}\left(\dfrac{F(y_2^2\tau_c|x)-F(y_2\tau_c|x)}{F(y_2\tau_c|x)-F(\tau_c|x)}\right),\\[.2cm]
p_{y_1,y_2,\tau_c}(x)&=&F(\tau_c|x)+\dfrac{F(\tau_c|x)-F(y_1\tau_c|x)}{y_1^{-1/\gamma_{y_2,\tau_c}(x)}-1},\\[.2cm]
F_{y_1,y_2,\tau_c}(t|x)&=&F(t\wedge\tau_c|x)+\left(p_{y_1,y_2,\tau_c}(x)-F(\tau_c|x)\right)\left(1-\left[\dfrac{t}{\tau_c}\vee 1\right]^{-1/\gamma_{y_2,\tau_c}(x)}\right).
\end{eqnarray*}

\section{Asymptotic properties} \label{section_asymptotics}

\subsection{The Beran estimator} \label{Beran}

In this subsection, we study the asymptotic properties of the \cite{Beran1981} estimator $F_n(t|x)$ for $t \le \tau_c(x)$, since most of the quantities involved in the previous section rely on increments of this estimator.  Previous asymptotic results for this estimator were restricted to $t \le \tau$ for some $\tau$ strictly less than $\tau_c(x)$ (see e.g.\ \cite{Gonzalez1994} or \cite{VanKeilegom1997}).  We extend these results to the full support of the censoring distribution thanks to an assumption on the censoring distribution (see assumption $(\mathcal{A})$ below).  For purposes of generality of exposition, in this subsection we allow $\tau_c(x)$ to depend on $x$. We also assume from now on that $x\in S_X$ defines a fixed reference position.

Due to the regression context, we need some H\"older-type conditions on the distribution functions $H$, $H^u$ and $F$ and on the density function $f$ of the covariate $X$. Let $\Vert \cdot \Vert$ be the Euclidean norm in $\mathbb{R}^p$.

\bigskip
\noindent
{\bf Assumption $(\mathcal{H})$.} There exist $0<\eta,\eta' \le 1$ and $c>0$ such that for any $t,s \in\mathbb{R}$ and any $x_1,x_2\in\mathbb{R}^p$,  \\
$(\mathcal{H}.1)$ $|f(x_1)-f(x_2)|\leq c\Vert x_1-x_2\Vert^\eta$, \\
$(\mathcal{H}.2)$ $|H(t|x_1)-H(t|x_2)|\leq c\Vert x_1-x_2\Vert^\eta$, \\
$(\mathcal{H}.3)$ $|H^u(t|x_1)-H^u(s|x_2)|\leq c(\Vert x_1-x_2\Vert^\eta+|t-s|^{\eta'})$, \\
$(\mathcal{H}.4)$ $|F(t|x)-F(s|x)|\leq c|t-s|^{\eta'}$.

\bigskip
Also, some common assumptions on the kernel function need to be imposed.

\bigskip
\noindent
{\bf Assumption $({\cal K})$.} Let $K$ be a bounded density function in $\mathbb R^p$ with support $S_K$ and suppose that $\int_{S_K}\Vert u\Vert K(u)du<\infty$. 

\bigskip
As a preliminary result, we show the rate of convergence of the following estimators of the functions $H(t|x)$ and $H^u(t|x)$ :
\begin{eqnarray*}
H_n(t|x) &=& \sum_{i=1}^n W_h(x-X_i) \ind_{\{T_i \le t\}}, \\
H_n^u(t|x) &=& \sum_{i=1}^n W_h(x-X_i) \ind_{\{T_i \le t,\delta_i=1\}}.
\end{eqnarray*}

\begin{lemma}
\label{lemma_rates}
Assume $(\mathcal{H}.1), (\mathcal{H}.2), (\mathcal{H}.3), ({\cal K})$  and $nh^{2\eta+p} |\log h|^{-1}=\mathcal{O}(1)$. Then, if $f(x)>0$,
\begin{eqnarray*}
\label{rate1}
\sup_{t\in\mathbb{R}}\left|H_n(t|x)-H(t|x)\right|=\mathcal{O}_\mathbb{P}\left((nh^p)^{-1/2}|\log h|^{1/2}\right),\\
\sup_{t\in\mathbb{R}}\left|H^u_n(t|x)-H^u(t|x)\right|=\mathcal{O}_\mathbb{P}\left((nh^p)^{-1/2}|\log h|^{1/2}\right). \nonumber
\end{eqnarray*}
\end{lemma}

We now derive an asymptotic representation for $F_n(t|x)-F(t|x)$ for $t \le \tau_c(x)$.  For this, we need one more assumption, that is crucial in order to obtain the representation upto (and including) the point $\tau_c(x)$.  

\bigskip
\noindent
{\bf Assumption $(\mathcal{A})$.} The distribution functions $F(\cdot|x)$ and $G(\cdot|x)$ are respectively continuous on $(-\infty,\tau_c(x)]$ and $(-\infty,\tau_c(x))$, and $
G(\tau_c(x)|x)-G(\tau_c(x)^-|x)>0$.

\medskip
\begin{remark}
\label{remark_H}
Assumption $(\mathcal{A})$ clearly implies that both the functions $H(t|x)$ and $H^u(t|x)$ are continuous for $t<\tau_c(x)$. Although $H(\tau_c(x)|x)-H(\tau_c(x)^-|x)>0$, one can show that $H^u(\cdot|x)$ is continuous at $\tau_c(x)$. Indeed, it is sufficient to see that under the model assumptions, $\mathbb{P}(Y=C|X=x)=0$, and hence
\begin{eqnarray*}
H^u(\tau_c(x)|x)=\mathbb{P}(Y\leq\tau_c(x),Y\leq C|X=x)&=&\mathbb{P}(Y\leq\tau_c(x),Y < C|X=x)\\
&=&\mathbb{P}(Y<\tau_c(x),Y < C|X=x)\\
&=&H^u(\tau_c(x)^-|x).
\end{eqnarray*}
\end{remark}

\begin{theorem}
\label{proposition_iid}
Assume $(\mathcal{A})$. Then, under the assumptions of Lemma \ref{lemma_rates} and for any $\tau<\tau_c(x)$, we have for $\tau\leq t \leq\tau_c(x)$,
\begin{eqnarray*}
F_n(t|x)-F(t|x)=\sum_{i=1}^n W_h(x-X_i) g(t,T_i,\delta_i|x) + r_n(t|x),
\end{eqnarray*}
where
\begin{eqnarray*}
g(t,T_i,\delta_i|x)&=&(1-F(t|x)) \left\{\int_{-\infty}^t\dfrac{\ind_{\{T_i < s\}}-H(s^-|x)}{(1-H(s^-|x))^2}dH^u(s|x)+\dfrac{\ind_{\{T_i \leq t,\delta_i=1\}}-H^u(t|x)}{1-H(t^-|x)}\right.\\
&&\left.-\int_{-\infty}^t\dfrac{\ind_{\{T_i \leq s,\delta_i=1\}}-H^u(s|x)}{(1-H(s^-|x))^2}dH(s^-|x)\right\},
\end{eqnarray*}
and
\begin{eqnarray*}
\sup_{\tau\leq t\leq\tau_c(x)}|r_n(t|x)|=\mathcal{O}_\mathbb{P}((nh^p)^{-3/4} |\log h|^{3/4}).
\end{eqnarray*}
\end{theorem}

This allows us to obtain the main result of this subsection, which is the weak convergence of the Beran estimator $F_n(t|x)$ as a process in $\ell^\infty[\tau,\tau_c(x)]$ for any $\tau < \tau_c(x)$ and for fixed $x$.  Here, for any set $S$, the space $\ell^\infty(S)$ is the space of bounded functions defined on $S$ endowed with the uniform norm. 

\begin{theorem}
\label{theorem_beran}
Assume $(\mathcal{A})$, $(\mathcal{H})$ and $(\mathcal{K})$, and assume that $f(x)>0$, $nh^p |\log h|^{-3} \to \infty$ and $nh^{2\eta+p-q} |\log h|^{-1} =\mathcal{O}(1)$ for some $q>0$. Then, for any $\tau < \tau_c(x)$, the process
\begin{eqnarray*}
\label{process}
\left\{(nh^p)^{1/2}(F_n(t|x)-F(t|x)),\quad t\in [\tau,\tau_c(x)]\right\},
\end{eqnarray*}
converges weakly in $\ell^\infty[\tau,\tau_c(x)]$ to a continuous mean-zero Gaussian process $Z(\cdot|x)$ with covariance function
\begin{eqnarray*}
\Gamma(t,s|x)=\dfrac{\Vert K\Vert_2^2}{f(x)}(1-F(t|x))(1-F(s|x))\int_{-\infty}^{t\wedge s}\dfrac{dH^u(y|x)}{(1-H(y^-|x))^2}.
\end{eqnarray*}
\end{theorem}
\vspace{.5cm}
\noindent

\subsection{The proposed estimators} \label{ours}

In this subsection, we develop the asymptotic properties of the estimators $\widehat\gamma(x)$, $\widehat p(x)$ and $\widehat F(t|x)$ defined in (\ref{gammahat}), (\ref{phat}) and (\ref{Fhat}). In order to derive their large sample properties, we need some additional assumptions on the model. First, it is important  to ensure that the right endpoint $\tau_c(x')$ does not depend on $x'$, and second that assumption $(\mathcal{A})$ is verified in a small neighbourhood of $x$.

\bigskip
\noindent
{\bf Assumption $(\mathcal{T})$.} The right endpoint $\tau_c(x')$ does not depend on the vector $x'$, and is henceforth denoted by $\tau_c$.

\bigskip 
\noindent
{\bf Assumption $(\mathcal{A}')$.} Let $B(x,r)$ be an open ball centered at $x\in\mathbb{R}^p$ of radius $r>0$ with respect to the norm $\Vert\cdot\Vert$. Then, there exists a $\delta>0$ such that $B(x,\delta)\subset S_X$ and $\inf_{x' \in B(x,\delta)} [G(\tau_c|x') - G(\tau_c^-|x')] >0$. 

\bigskip
\noindent
Both assumptions $(\mathcal{A}')$ and $(\mathcal{T})$ are needed for the estimation of $\tau_c$ and ensure that the speed of convergence of $\tau_n$ is fast enough for our method. 

\begin{theorem}
\label{proposition_gamma}
Assume $(\mathcal{A}')$ and $(\mathcal{T})$. Under the conditions of Theorem \ref{theorem_beran} and for any $y_2\in(0,1)$ such that 
\begin{eqnarray*}
b(y_2|x)=\dfrac{F(y_2^ 2\tau_c|x)-F(y_2\tau_c|x)}{F(y_2\tau_c|x)-F(\tau_c|x)}\neq 0,
\end{eqnarray*}
we have 
\begin{eqnarray*}
(nh^p)^{1/2}\big(\widehat\gamma(x)-\gamma_{y_2,\tau_c}(x)\big)\overset{d}{\longrightarrow}\mathcal{N}\big(0,\sigma^2_{\gamma,\tau_c}(x)\big),
\end{eqnarray*}
with 
\begin{eqnarray*}
\sigma^2_{\gamma,\tau_c}(x)=\left[\dfrac{\phi_{y_2}(b(y_2|x))}{a(y_2|x)}\right]^2\sum_{i,j=0}^2c_{i,j}\Gamma(y_2^i\tau_c,y_2^j\tau_c|x),
\end{eqnarray*}
where $a(y_2|x)=F(y_2\tau_c|x)-F(\tau_c|x)$, $\phi_{y_2}(x)=\dfrac{\log(y_2)}{x\log(x)^2}$, $c_{i,j}=c_{j,i}$ and \\ 

\begin{center}
$\begin{array}{ll}
\quad c_{0,0}=b(y_2|x)^2,  \quad \quad \quad \quad &c_{0,1}=-b(y_2|x)(1+b(y_2|x)), \\
\quad c_{1,1}=[1+b(y_2|x)]^2, &c_{0,2}=b(y_2|x),\\
\quad c_{1,2}=-(1+b(y_2|x)), &c_{2,2}=1. \\
\end{array}$
\end{center}
\end{theorem}

\bigskip

\begin{theorem}
\label{theorem_p}
Under the assumptions of Theorem \ref{proposition_gamma} and for any $y_1 \in (0,1)$ such that 
\begin{eqnarray*}
y_1^{-1/\gamma_{y_2,\tau_c}(x)}-1\neq 0,
\end{eqnarray*}
we have 
\begin{eqnarray*}
(nh^p)^{1/2}\big(\widehat p(x)-p_{y_1,y_2,\tau_c}(x)\big)\overset{d}{\longrightarrow}\mathcal{N}\big(0,\sigma^2_{p,\tau_c}(x)\big),
\end{eqnarray*}
with
\begin{eqnarray*}
\sigma^2_{p,\tau_c}(x)=\sum_{i,j=0}^3d_ig_{i,j}d_j,
\end{eqnarray*}
where
\begin{eqnarray*}
d_0&=&\dfrac{y_1^{-1/\gamma_{y_2,\tau_c}(x)}}{y_1^{-1/\gamma_{y_2,\tau_c}(x)}-1}-\dfrac{\log(y_1)y_1^{-1/\gamma_{y_2,\tau_c}(x)}}{(\gamma_{y_2,\tau_c}(x)(y_1^{-1/\gamma_{y_2,\tau_c}(x)}-1))^2}(F(\tau_c|x)-F(y_1\tau_c|x))\phi_2(b(y_2|x)),\\
d_1&=&-\dfrac{1}{y_1^{-1/\gamma_{y_2,\tau_c}(x)}-1},\\
d_2&=&\dfrac{\log(y_1)y_1^{-1/\gamma_{y_2,\tau_c}(x)}}{(\gamma_{y_2,\tau_c}(x)(y_1^{-1/\gamma_{y_2,\tau_c}(x)}-1))^2}\dfrac{F(\tau_c|x)-F(y_1\tau_c|x)}{F(\tau_c|x)-F(y_2\tau_c|x)}\phi_{y_2}(b(y_2|x)) \\
&& \times \left(1-\dfrac{1}{F(y_2\tau_c|x)-F(\tau_c|x)}\right),\\
d_3&=&\dfrac{\log(y_1)y_1^{-1/\gamma_{y_2,\tau_c}(x)}}{(\gamma_{y_2,\tau_c}(x)(y_1^{-1/\gamma_{y_2,\tau_c}(x)}-1))^2}\dfrac{F(\tau_c|x)-F(y_1\tau_c|x)}{F(\tau_c|x)-F(y_2\tau_c|x)}\phi_{y_2}(b(y_2|x)),
\end{eqnarray*}
and $g_{i,j}=g_{j,i}$ with
\begin{eqnarray*}
g_{i,j}=\left\{\begin{array}{ll}
\Gamma(y_1^i\tau_c,y_1^j\tau_c|x),&\text{if}\quad i\leq 1,\quad j\leq 1,\\
\Gamma(y_1^i\tau_c,y_2^{j-1}\tau_c|x),&\text{if}\quad i\leq 1,\quad 2\leq j\leq 3,\\
\Gamma(y_2^{i-1}\tau_c,y_2^{j-1}\tau_c|x),&\text{if}\quad 2\leq i\leq 3,\quad 2\leq j\leq 3.
\end{array}
\right.
\end{eqnarray*}
\end{theorem}

\bigskip

\begin{theorem}
\label{theorem_F}
Under the assumptions of Theorem \ref{theorem_p} and for any $a \in \mathbb{R}$, the process 
\begin{eqnarray*}
\left\{(nh^p)^{1/2}(\widehat{F}(t|x)-F_{y_1,y_2,\tau_c}(t|x)),\quad t\in[a,+\infty)\right\}
\end{eqnarray*}
converges weakly in $\ell^\infty[a,\infty)$ to a continuous mean-zero Gaussian process $\widetilde{Z}(\cdot|x)$ with covariance function
\begin{eqnarray*}
\Delta_{y_1,y_2,\tau_c}(t,s|x)= \sum_{i,j=0}^3e_i(t)g_{i,j}e_j(s)+\sum_{i=0}^3 \big(\widetilde{g}_i(t)e_i(s)+\widetilde{g}_i(s)e_i(t)\big),
\end{eqnarray*}
where for any $t\in\mathbb{R}$,
\begin{eqnarray*}
&&\widetilde{g}_0(t)=\Gamma(\tau_c,t\wedge\tau_c|x),\quad\hspace{.35cm} \widetilde{g}_1(t)=\Gamma(y_1\tau_c,t\wedge\tau_c|x),\\
&&\widetilde{g}_2(t)=\Gamma(y_2\tau_c,t\wedge\tau_c|x),\quad\widetilde{g}_3(t)=\Gamma(y_2^2\tau_c,t\wedge\tau_c|x),
\end{eqnarray*}
and 
\begin{eqnarray*}
e_0(t)&=&-\dfrac{p-F(\tau_c|x)}{\gamma^2_{y_2,\tau_c}(x)}\log\left(\frac{t}{\tau_c}\vee 1\right)\left[\frac{t}{\tau_c}\vee 1\right]^{-1/\gamma_{y_2,\tau_c}(x)}\phi_{y_2}\left(b(y_2|x)\right)\dfrac{b(y_2|x)}{a(y_2|x)}\\
&&+\left(1-\left[\frac{t}{\tau_c}\vee 1\right]^{-1/\gamma_{y_2,\tau_c}(x)}\right)\left(\dfrac{1}{y^{-1/\gamma_{y_2,\tau_c}(x)}-1} \right.\\[.1cm]
&&  \quad\quad  -c(y_1,y_2|x)\psi_{y_1}(\gamma_{y_2,\tau_c}(x))\phi_{y_2}(b(y_2|x))b(y_2|x)\Big),\\
e_1(t)&=&\dfrac{1-\left[\frac{t}{\tau_c}\vee 1\right]^{-1/\gamma_{y_2,\tau_c}(x)}}{y^{-1/\gamma_{y_2}}-1},\\
e_2(t)&=&\dfrac{p-F(\tau_c|x)}{\gamma^2_{y_2,\tau_c}(x)}\log\left(\frac{t}{\tau_c}\vee 1\right)\left[\frac{t}{\tau_c}\vee 1\right]^{-1/\gamma_{y_2,\tau_c}(x)}\phi_{y_2}\left(b(y_2|x)\right)\dfrac{1+b(y_2|x)}{a(y_2|x)}\\
&&+\left(1-\left[\frac{t}{\tau_c}\vee 1\right]^{-1/\gamma_{y_2,\tau_c}(x)}\right)c(y_1,y_2|x)\psi_{y_1}(\gamma_{y_2,\tau_c}(x))\phi_{y_2}(b(y_2|x))(1-b(y_2|x)),\\
e_3(t)&=&\dfrac{p-F(\tau_c|x)}{\gamma^2_{y_2,\tau_c}(x)}\log\left(\frac{t}{\tau_c}\vee 1\right)\left[\frac{t}{\tau_c}\vee 1\right]^{-1/\gamma_{y_2,\tau_c}(x)}\phi_{y_2}\left(b(y_2|x)\right)\dfrac{1}{a(y_2|x)}\\
&&+\left(1-\left[\frac{t}{\tau_c}\vee 1\right]^{-1/\gamma_{y_2,\tau_c}(x)}\right)c(y_1,y_2|x)\psi_{y_1}(\gamma_{y_2,\tau_c}(x))\phi_{y_2}(b(y_2|x)),
\end{eqnarray*}
with $c(y_1,y_2|x)=\dfrac{F(y_1\tau_c|x)-F(\tau_c|x)}{F(y_2\tau_c|x)-F(\tau_c|x)}$.
\end{theorem}

\section{Simulations}\label{section_simulation}

In this section we study the finite sample performance of our estimators by means of a simulation study.  We assume throughout that the covariate $X$ is uniformly distributed on the interval $[0,1]$.  For the non-cured subjects ($Y<\infty$), we consider three models for the distribution $F_0(\cdot|x)$ for $X=x$ : a generalized extreme value distribution (GEV), a generalized Pareto distribution (GPD) and a Fr\'echet distribution, sharing the same conditional extreme value index $\gamma(x)$ given by $\gamma(x)=(x+1)/2$.  The conditional probability of not being cured $p(x)$ follows a logistic model given by 
\begin{eqnarray*}
p(x)=\dfrac{\exp(\beta_1+\beta_2(2x-1))}{1+\exp(\beta_1+\beta_2(2x-1))},
\end{eqnarray*}  
with $(\beta_1,\beta_2)=(0.4,2)$.   The censoring time $C$ is independent of $X$ and independent of $Y$ given $X$, and it is uniformly distributed on the interval $[0,\tau_c]$ with probability $1-\varepsilon>0$, and fixed to $\tau_c$ otherwise.  For the right endpoint $\tau_c$, we consider a range of values between $\tau_{0.25}$ and $\tau_{0.95}$ : $\tau_{c,s}=\tau_{0.25}+s(\tau_{0.95}-\tau_{0.25})$, where $s\in[0,1]$ and $\tau_{\alpha}$ is the quantile of level $\alpha\in(0,1)$ for $F_0(.|x)$. 
\begin{figure}[H]
\centering
\includegraphics[scale=0.45]{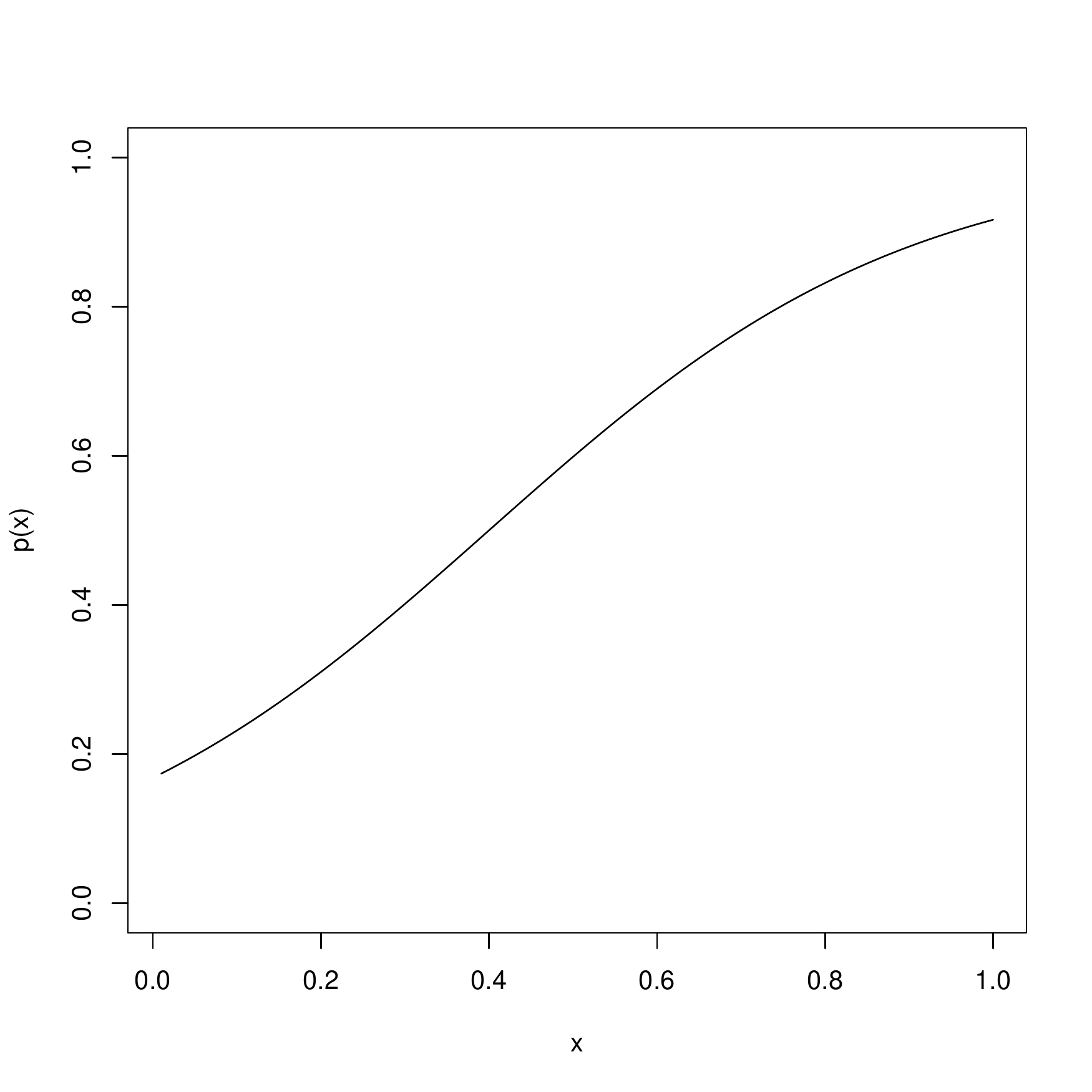}
\caption{Conditional probability of not being cured as a function of the covariate. }
\label{figure_p}
\end{figure}

Our simulations are based on datasets of size $n=2000$ and the procedure is repeated $N=100$ times.   To compute our estimators $\widehat \gamma(x)$, $\widehat p(x)$ and $\widehat F(t|x)$, we need to select an appropriate bandwidth $h_n$, a kernel function $K$ and tuning parameters $y_1$ and $y_2$.  The bandwidth $h_n$ is chosen by means of the function \texttt{dpik} from the \texttt{R}-package \texttt{KernSmooth}, and the kernel $K$ is the Epanechnikov kernel $K(u)=(3/4)(1-u^2)\ind_{\{|u|<1\}}$.   Note that by construction, $\widehat \gamma(x)$ might take negative values. In order to satisfy the constraint that $\gamma(x)>0$, we truncate the estimator from below by considering $\widehat \gamma(x) \vee 0.1$.  Next, the pair $(y_1,y_2)$ is selected in a data-driven way by means of a quadratic errors criterion between the possible pairs, namely
\begin{eqnarray*}
(y_1,y_2)=\underset{(z_1,z_2)\in\mathcal{G}^2}{\text{argmin}}\sum_{(z'_1,z'_2)\in\mathcal{G}^2}\left(\widehat p_{z_1,z_2}(x)-\widehat p_{z'_1,z'_2}(x)\right)^2,
\end{eqnarray*}
where $\mathcal{G}=\{0.25,0.27,\ldots,0.89\}$ is a grid of values, and $\widehat p_{z_1,z_2}(x)$ refers to the estimator $\widehat p(x)$ using the tuning parameters $z_1$ and $z_2$.

\begin{figure}[H]
\centering
\includegraphics[scale=0.55,page=3]{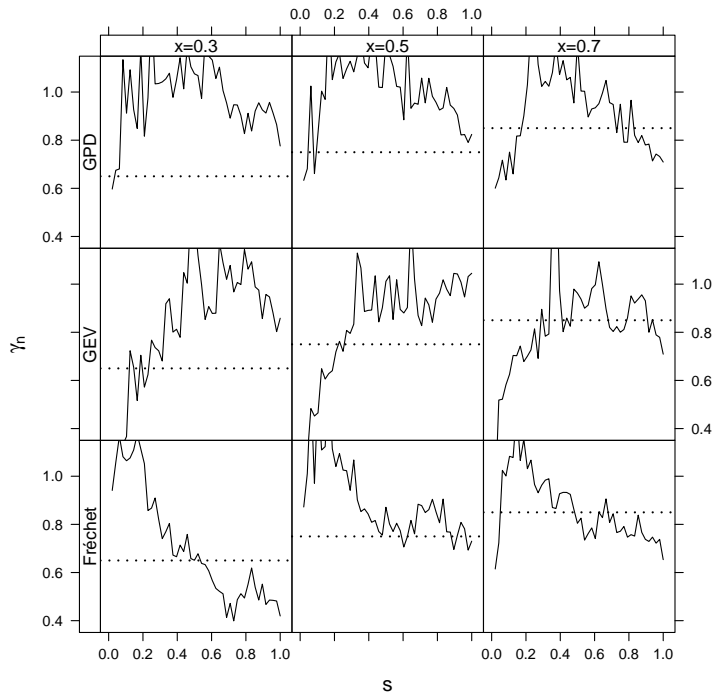}
\caption{Median results for $\widehat\gamma(x)$ for $x=0.3,0.5,0.7$ and for $\tau_c=\tau_{c,s}$ with $0 \le s \le 1$. The dotted horizontal line represents the true value of $\gamma(x)$.}
\label{figure_evi}
\end{figure}

We first study the performance of the estimator $\widehat \gamma(x)$ for $x=0.3, 0.5, 0.7$.  
From Figure \ref{figure_evi}, which shows the median of the estimator $\widehat\gamma(x)$,  we observe that the estimator is quite unstable and biased. For $x=0.3$, the estimator has difficulties approximating the true value of $\gamma(x)$. This can probably be explained by the high value of the cure rate ($1-p(x)\simeq 0.6$), which implies that the censoring proportion in the model is at least $60\%$.  Also note that the tuning parameter $y_2$ is chosen so as to optimize the estimation of $p(x)$, but the selected value of $y_2$ might not be the best value for the estimation of $\gamma(x)$. For $x=0.5$ and $x=0.7$, the smaller censoring percentage improves the accuracy of the estimator for all distributions. Among them, the GEV distribution shows the best results with a reduced bias, although all models suffer from a certain variability as a function of $s$. Nevertheless, we will see in the sequel that the variability of $\widehat\gamma(x)$ does not affect the accuracy of the estimators $\widehat p(x)$ and $\widehat F(t|x)$. 

\begin{figure}[H]
\centering
\includegraphics[scale=0.55,page=1,trim = 5cm 0cm 0cm 0cm]{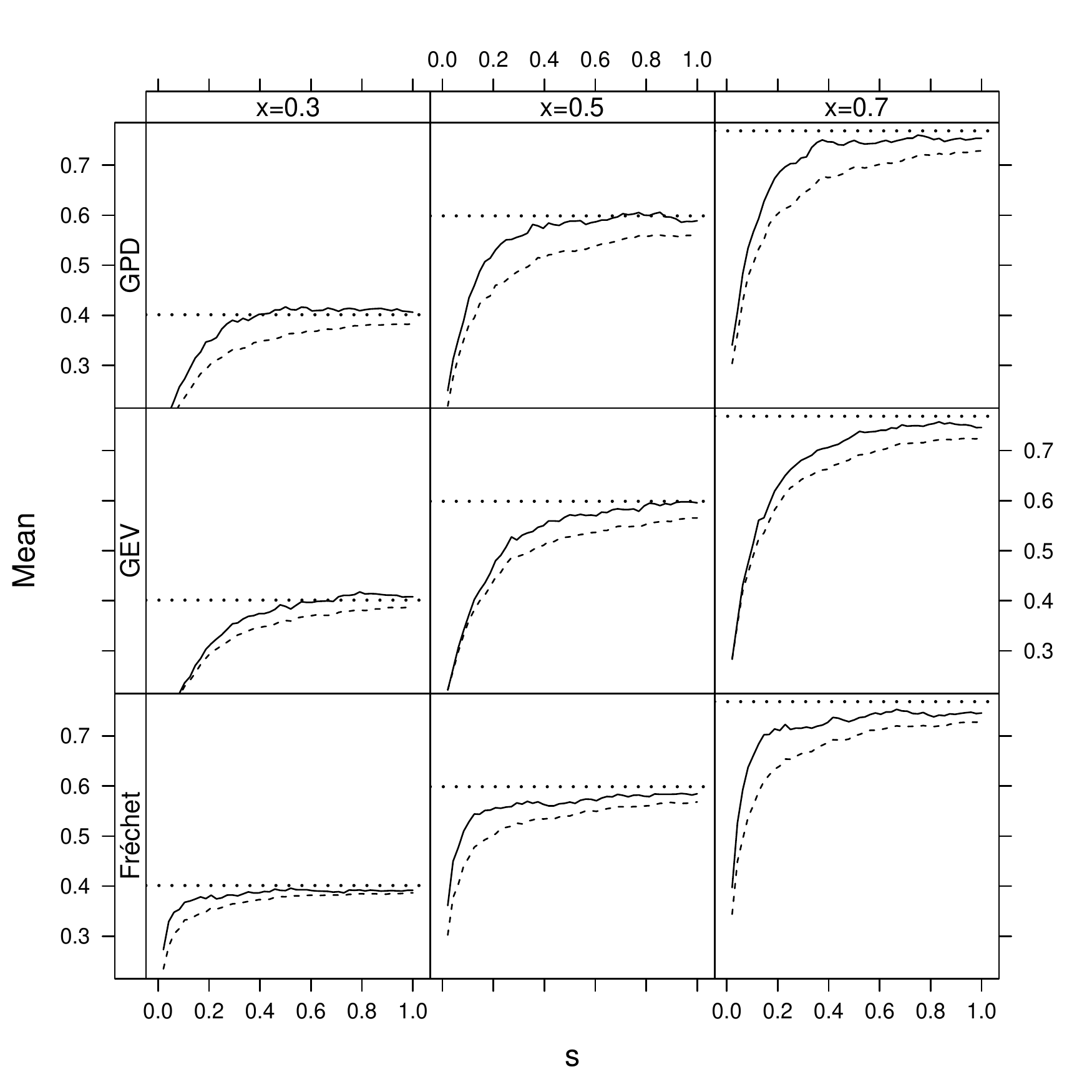}
\includegraphics[scale=0.55,page=2,trim = 0cm 0cm 5cm 0cm]{{datadriven_trellis}.pdf}
\caption{Mean (left) and MSE (right) of the Beran estimator $p_n(x)$ (dashed) and of our estimator $\widehat p(x)$ (solid) for $x=0.3,0.5,0.7$ and for $\tau_c=\tau_{c,s}$ with $0 \le s \le 1$. The dotted horizontal line represents the true value of $p(x)$.}
\label{figure_pn}
\end{figure}

Next, we compare our estimator $\widehat p(x)$ of the conditional non-cure rate with the Beran estimator $p_n(x) = F_n(\tau_n|x)$ as a function of $\tau_{c,s}$ with $0 \le s \le 1$.   To compare the two estimators we compute the mean and the mean squared error (MSE), which is given by 
\begin{eqnarray*}
\text{MSE}(x) = \dfrac{1}{N}\sum_{i=1}^N\left(\widehat p^{(i)}(x)-p(x)\right)^2.
\end{eqnarray*}
Here, $\widehat p^{(i)}(x)$ is our estimator $\widehat p(x)$ obtained with the $i$-th sample,  and similarly for the Beran estimator when $\widehat p^{(i)}(x)$ is replaced by $p^{(i)}_n(x)$.
Figure \ref{figure_pn} shows the mean and the MSE of the Beran estimator $p_n(x)$ and of our estimator $\widehat p(x)$. As it is theoretically expected, the Beran estimator always underestimates the true $p(x)$. On the contrary, $\widehat p(x)$ gets close to $p(x)$ for a wide range of $s$-values. Satisfactory results are indeed observed for $s\geq a$ with $a=0.2$ for the Fr\'echet distribution, $a=0.3$ for the GPD and $a=0.4$ for the GEV distribution. Our estimator particularly shows smooth curves that tend to stabilize around the true value for $s\geq a$. In terms of MSE, we observe in most cases lower or similar MSE values than for the Beran estimator, but for $x=0.3$ the estimators again suffer from a high censoring proportion compared to $x=0.5$ and $x=0.7$. Overall, the bias is reduced for $\widehat p(x)$ and the MSE curves show that the balance between bias correction and variability is still more advantageous for our approach. 

\begin{figure}[H]
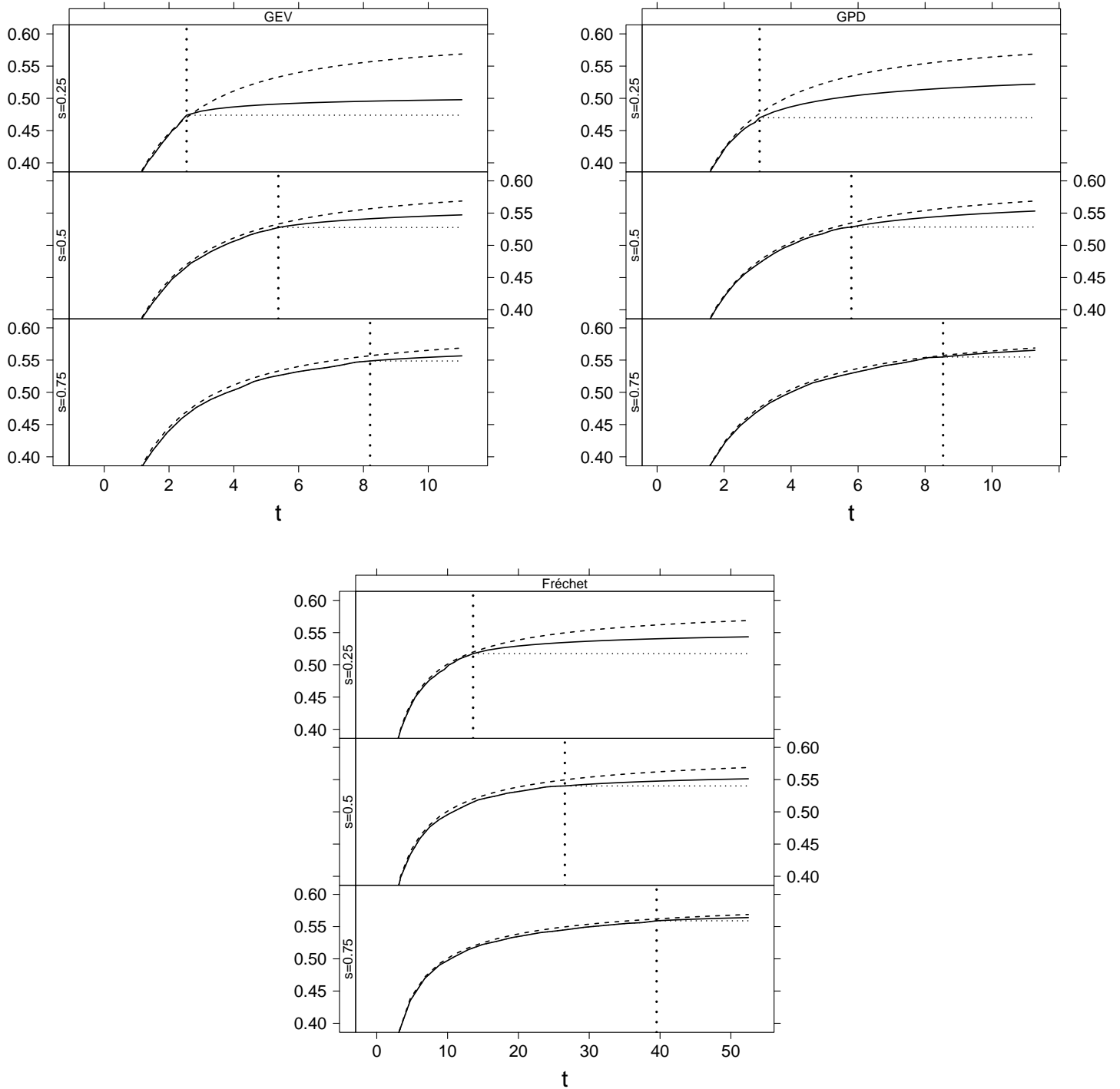

\centering
\includegraphics[scale=0.55,page=7,trim=5cm 0cm 0cm 1cm]{{datadriven_trellis}.pdf}
\includegraphics[scale=0.55,page=8,trim=0cm 0cm 5cm 1cm]{{datadriven_trellis}.pdf}
\includegraphics[scale=0.55,page=9,trim= 0cm 0cm 0cm 0cm]{{datadriven_trellis}.pdf}
\caption{Mean of $\widehat{F}(t|x)$ (solid curve) and $F_n(t|x)$ (dotted curve) for $x=0.5$, $t\in[\tau_{0.25},\tau_{0.95}]$ and $s=0.25,0.5,0.75$. The true distribution function $F(t|x)$ is given by the dashed  curve while the value of $\tau_c$ is shown by the vertical dotted line.}
\label{figure_df_0.5}
\end{figure}

Finally, Figure \ref{figure_df_0.5} shows the mean of our estimator $\widehat F(t|x)$ for $t\in[\tau_{0.25},\tau_{0.95}]$, where the vertical dotted line refers to the value of $\tau_c$. To the left of this line, $\widehat F(t|x)$ equals the Beran estimator $F_n(t|x)$, whereas to the right, the Beran estimator remains constant whereas our estimator makes use of the extrapolation approach to remain close to the true distribution.  Also note that the true distribution $F(t|x)$ grows to $1-p(x) \sim 0.6$ when $t$ tends to infinity, our estimator stays relatively close to this true curve (especially when $s$ increases), whereas the Beran estimator remains constant for $t \ge \tau_c$ at a value between 0.35 and 0.60 depending on the value of $s$ and on the underlying distribution.  

\section{Application} \label{section_data}

\begin{figure}[H]
\centering
\includegraphics[scale=0.4]{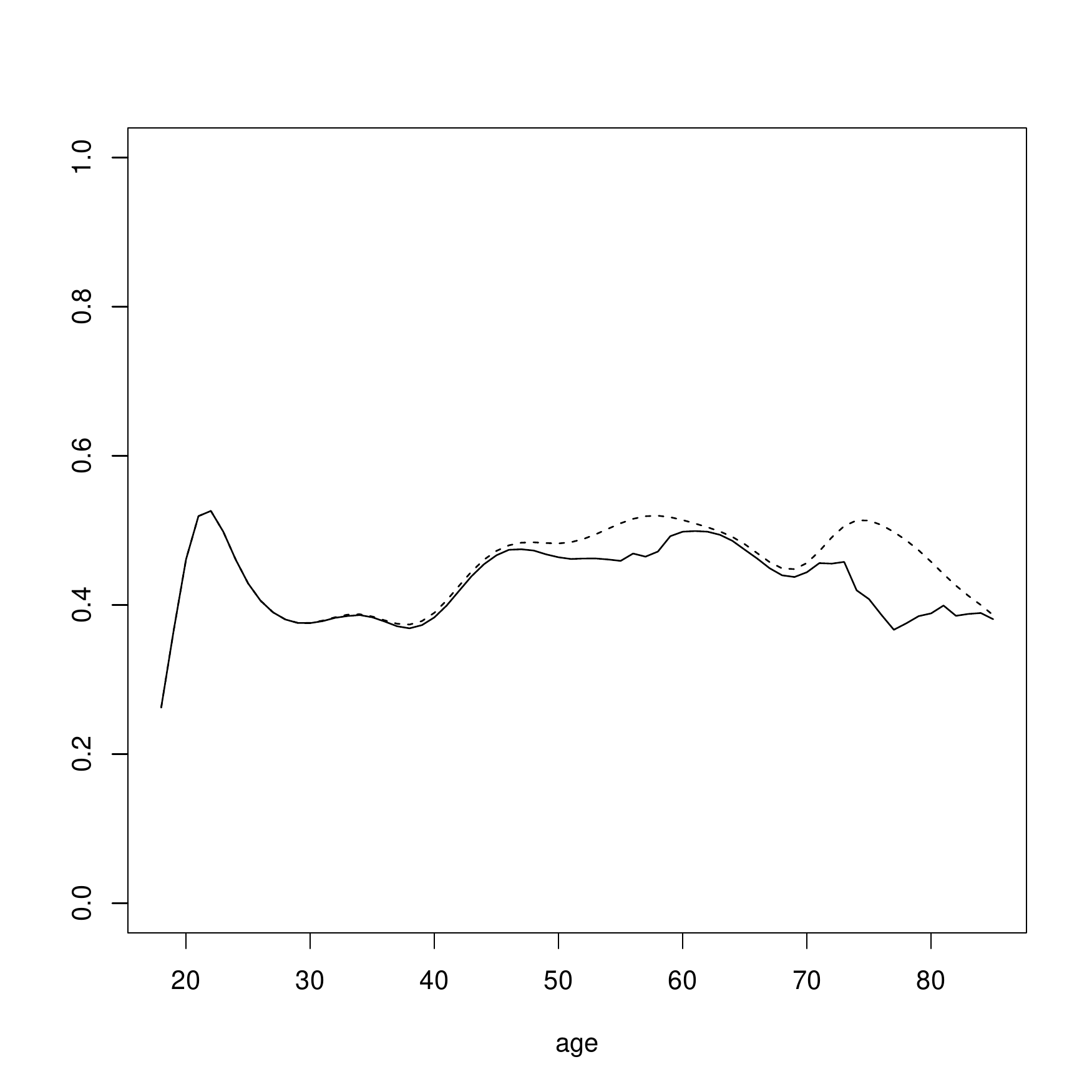}
\includegraphics[scale=0.4]{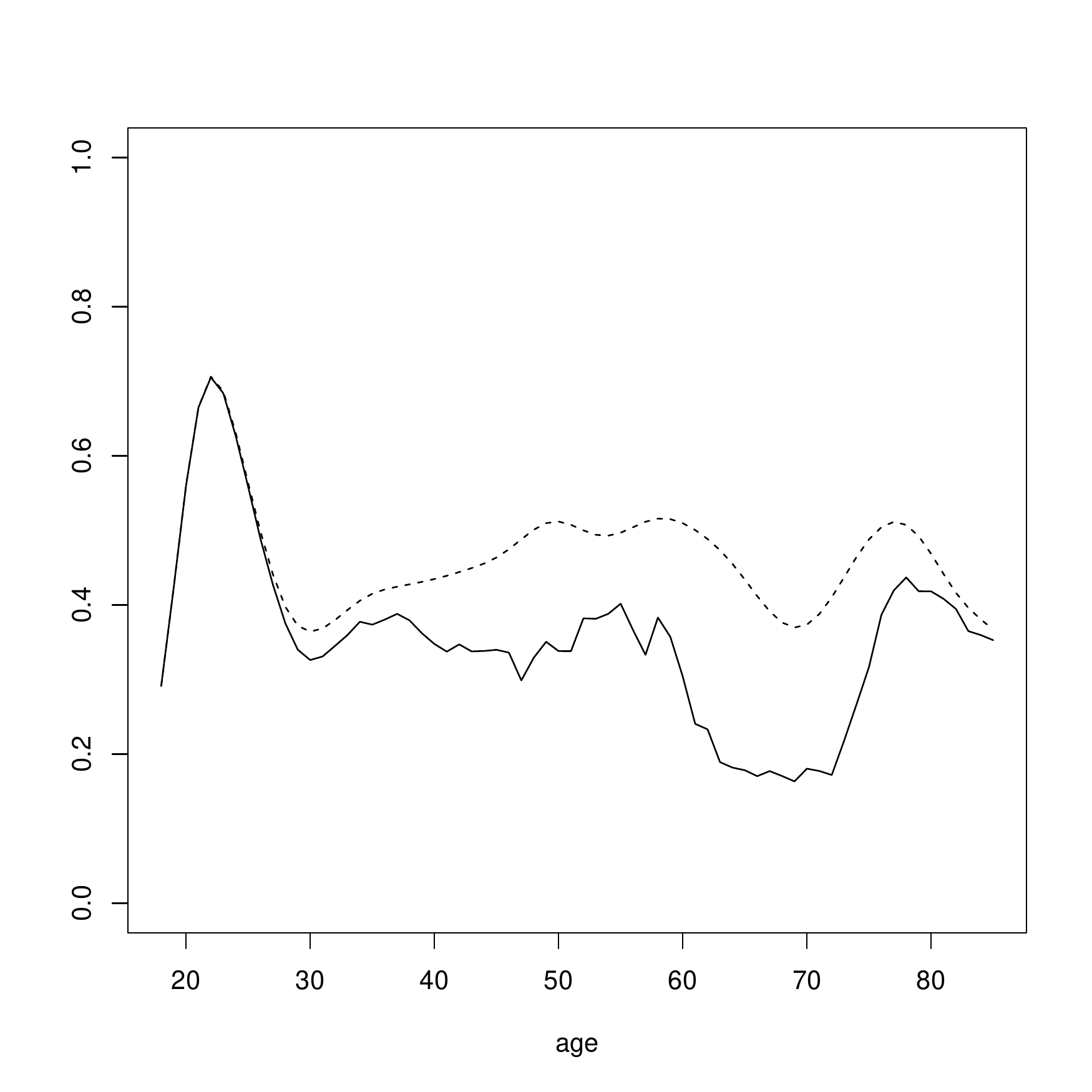}
\caption{Estimation of the cure rate $1-p(x)$ as a function of age by means of the Beran estimator $1-p_n(x)$ (dashed curve) and our estimator $1-\widehat p(x)$ (solid curve).  The event of interest is the recurrence of the cancer (left panel) and the death of the patients (right panel).}
\label{figure_colon}
\end{figure}

In this section we conduct a small study on the survival of colon cancer patients. The dataset is available in the \texttt{R}-package \texttt{survival} (called \texttt{colon}).  Patients in this study receive chemotherapy by means of either Levamisole, or 5-FU in addition to Levamisole.  Two events are of interest : the recurrence of the colon cancer and the death of the patients.  The sample consists of 929 patients, and for each patient the two event times, possibly right censored, are recorded.   Figure \ref{figure_colon}  shows the estimated cure rates for the two event types as a function of age.   Both the Beran estimator $1-p_n(x)$ and our estimator $1-\widehat p(x)$ are calculated.   

These figures clearly show two opposite phenomena. When the event of interest is the recurrence of the cancer, the discrepancy between the Beran estimator $1-p_n(x)$ and our estimator $1-\widehat p(x)$ is quite small, suggesting that the follow-up for this event type is sufficient.  Hence, recurrence of the cancer happens relatively soon, so that the recurrence rate can be estimated consistently with the naive Beran estimator, without the need to use extrapolation techniques. On the other hand, the discrepancy between the two estimators is much bigger when the event of interest is the death of the patients.  This strongly suggests that the follow-up is not sufficient in this case, meaning that a certain proportion of the patients die after the end of the study, and so extrapolation techniques are needed in order to correctly estimate the cure rate.

\section{Appendix: Proofs of the main results}

Before giving the proofs of the main results, we need to introduce some more notations mostly borrowed from the theory of weak convergence of empirical processes. For any class $\mathcal{F}$ of bounded and measurable functions over a metric space $(\mathcal{T},d)$ and any probability measure $Q$ and $\epsilon>0$, define the covering number $N(\mathcal{F},L_2(Q),\epsilon)$ as the minimal number of $L_2(Q)-$balls of radius $\epsilon$ needed to cover $\mathcal{F}$. We say that the class $\mathcal{F}$ is \textit{VC} if one can find $A>0$ and $\nu>0$ such that for any probability measure $Q$ and $\epsilon>0$,
\begin{eqnarray*}
N(\mathcal{F},L_2(Q),\epsilon\Vert F\Vert_{Q,2})\leq\left(\dfrac{A}{\epsilon}\right)^\nu,
\end{eqnarray*} 
where $0<\Vert F\Vert^2_{Q,2}=\int F^2 dQ<\infty$ and $F$ is an envelope function of the class ${\cal F}$. Additionally, we also define the uniform entropy integral as
\begin{eqnarray*}
J(\delta,\mathcal{F},L_2) =\int_0^\delta\sqrt{\log\sup_{\mathcal{Q}}N(\mathcal{F},L_2(Q),t\Vert F\Vert_{Q,2})}\, dt,
\end{eqnarray*} 
where $\mathcal{Q}$ is the set of all probability measures $Q$.

\subsection{Proofs of Section \ref{Beran}}

\noindent
{\bf Proof of Lemma \ref{lemma_rates}.}  The arguments for $H^u_n(\cdot|x)$ and $H_n(\cdot|x)$ are exactly the same and thus we only give the proof for the latter. Let us first show that
\begin{eqnarray}
\label{proof_rate1}
\sup_{t\in\mathbb{R}}\left|\mathbb{E}\left[K_h(x-X)\ind_{T\leq t}\right]-f(x)H(t|x)\right|=\mathcal{O}(h^\eta).
\end{eqnarray} 
Indeed, for any $t$ and $n$ large enough,
\begin{eqnarray*}
&& \mathbb{E}\left[K_h(x-X)\ind_{T\leq t}\right] \\
&& = \int_{S_X}h^{-p}K((x-u)/h)H(t|u)f(u)du\\
&& =\int_{x-hS_X}K(u)H(t|x-hu)f(x-hu)du\\
&& =H(t|x)f(x)+\int_{x-hS_X}K(u)[H(t|x-hu)f(x-hu)-H(t|x)f(x)]du,
\end{eqnarray*}
yielding
\begin{eqnarray*}
&&\left|\mathbb{E}\left[K_h(x-X)\ind_{T\leq t}\right]-f(x)H(t|x)\right|\\
&&\leq \int_{S_K}K(u)\left|H(t|x-hu)f(x-hu)-H(t|x)f(x)\right|du\\
&&\leq \int_{S_K}K(u)\left|f(x-hu)-f(x)\right|du+f(x)\int_{S_K}K(u)\left|H(t|x-hu)-H(t|x)\right|du\\
&&\leq h^\eta c\int_{S_K}K(u)\Vert u\Vert^\eta du(1+f(x)),
\end{eqnarray*}
and thus (\ref{proof_rate1}) is proven. 

Clearly, the class $\mathcal{G}=\left\{u\in\mathbb{R}\to\ind_{\{u\leq t\}},\,t\in\mathbb{R}\right\}$ is \textit{VC} since the class of subsets \linebreak $\left\{(u,t)\in\mathbb{R}^2,\,u\leq t\right\}$ forms a \textit{VC} class of sets (see \citealt{Vaart1996}, for more details about these concepts), meaning that $\mathcal{G}$ is \textit{VC} and there exist $A,\nu>0$ such that for any $\epsilon>0$, 
\begin{eqnarray*}
N(\mathcal{G},L_2(Q),\epsilon)\leq\left(\dfrac{A}{\epsilon}\right)^\nu.
\end{eqnarray*}
In particular, for $\mathcal{G}_n=\left\{(u,w)\in\mathbb{R}\times\mathbb{R}^p\to K_h(x-w)\ind_{\{u\leq t\}},\,t\in\mathbb{R}\right\}$ we have that 
\begin{eqnarray*}
N(\mathcal{G}_n,L_2(Q),\epsilon\Vert K_h(x-.)\Vert_{Q,2})\leq\left(\dfrac{A}{\epsilon}\right)^\nu,
\end{eqnarray*}
since we only update the previous sets with one single function and only one ball is needed to cover the class $\left\{w\in\mathbb{R}^p\to K_h(x-w)\right\}$ whatever the measure $Q$. From the concentration inequality (2.2) in Theorem 2.1 in \cite{Gine2002}, it follows that for $U=\Vert K\Vert_\infty = \sup_u K(u)$, for $\sigma^2=\sigma^2_n=h^{-p}\Vert K\Vert_\infty^2$ and for some universal constant $B>0$, we have for $\varepsilon_n=(nh^{-p} |\log h|)^{1/2}$,
\begin{eqnarray*}
\label{proof_rate3}
\nonumber&&\mathbb{E}\left[\sup_{t\in\mathbb{R}}\left|\sum_{i=1}^n\left( K_h(x-X_i)\ind_{\{T_i\leq t\}} - n\mathbb{E}\left[K_h(x-X)\ind_{\{T\leq t\}}\right]\right)\right|\right]\\
&&\leq B\left[\nu U\log(AU/\sigma)+\left(\nu n\sigma^2\log(AU/\sigma)\right)^{1/2}\right]=\mathcal{O}(\varepsilon_n).
\end{eqnarray*}
In particular, this implies that
\begin{eqnarray*}
\sup_{t\in\mathbb{R}}\left|\dfrac{1}{n}\sum_{i=1}^n K_h(x-X_i)\ind_{\{T_i\leq t\}}-f(x)H(t|x)\right|=\mathcal{O}_\mathbb{P}(n^{-1}\varepsilon_n+h^\eta)=\mathcal{O}_\mathbb{P}((nh^p)^{-1/2}|\log h|^{1/2}).
\end{eqnarray*}
The same arguments can be used to show that
\begin{eqnarray}
\label{kerneldensity}
\left|\dfrac{1}{n}\sum_{i=1}^n K_h(x-X_i)-f(x)\right|=\mathcal{O}_\mathbb{P}((nh^p)^{-1/2} |\log h|^{1/2}),
\end{eqnarray}
which shows the result.\hfill$\Box$

\bigskip

\noindent
{\bf Proof of Theorem \ref{proposition_iid}.} The proof mainly follows the same ideas as in the proof of Theorem 3.1 in \cite{VanKeilegom1997}, but we extend their result to $t\leq\tau_c(x)$ by making use of the singularity of $G(\cdot|x)$ at $\tau_c(x)$. To do so, define
\begin{eqnarray*}
G_n(t|x)=1-\prod_{T_{(i)}\leq t}\left(1-\dfrac{W_h(x-X_{(i)})}{1-\sum_{j=1}^{i-1}W_h(x-X_{(j)})}\right)^{1-\delta_{(i)}},
\end{eqnarray*}
which is nothing but the Beran estimator for the censoring distribution. It turns out that
\begin{eqnarray*}
(1-F_n(t|x))(1-G_n(t|x))&=&\prod_{T_{(i)}\leq t}\left(1-\dfrac{W_h(x-X_{(i)})}{1-\sum_{j=1}^{i-1}W_h(x-X_{(j)})}\right)\\
&=&\dfrac{\prod_{T_{(i)}\leq t}\left(1-\sum_{j=1}^{i}W_h(x-X_{(j)})\right)}{\prod_{T_{(i)}\leq t}\left(1-\sum_{j=1}^{i-1}W_h(x-X_{(j)})\right)}\\
&=&1-H_n(t|x).
\end{eqnarray*}
Next, the equality $dH^u_n(\cdot|x)=(1-G_n(\cdot^-|x))dF_n(\cdot|x)$ gives us
\begin{eqnarray*}
\Lambda_n(t|x) &=& \int_{-\infty}^t\dfrac{dF_n(s|x)}{1-F_n(s^-|x)} \\
&=&\int_{-\infty}^t\dfrac{F_n(s^-|x)-F_n(s|x)}{(1-F_n(s^-|x))(1-F_n(s|x))}dF_n(s|x)-\log(1-F_n(t|x))\\
&=&\mathcal{O}_\mathbb{P}((nh^p)^{-1})-\log(1-F_n(t|x)),
\end{eqnarray*} 
since $\sup_{t\in\mathbb{R}}|F_n(t^-|x)-F_n(t|x)|=\mathcal{O}_\mathbb{P}((nh^p)^{-1})$. By continuity, we have $1-F(t|x)=\exp\left(-\Lambda(t|x)\right)$, which yields the identity
\begin{eqnarray*}
F_n(t|x)-F(t|x)=\left[e^{-\Lambda(t|x)}-e^{-\Lambda_n(t|x)}\right]+\mathcal{O}_\mathbb{P}((nh^p)^{-1}).
\end{eqnarray*}
Using a Taylor expansion of second order, we obtain
\begin{eqnarray*}
F_n(t|x)-F(t|x)=(1-F(t|x))(\Lambda_n(t|x)-\Lambda(t|x))+R_{n,1}(t|x)+\mathcal{O}_\mathbb{P}((nh^p)^{-1}),
\end{eqnarray*}
where
\begin{eqnarray*}
R_{n,1}(t|x)&=&-\dfrac{1}{2}\exp(\widetilde{\Lambda}(t|x))\left[\Lambda_n(t|x)-\Lambda(t|x)\right]^2,
\end{eqnarray*}
with $\widetilde{\Lambda}(t|x)$ between $\Lambda_n(t|x)$ and $\Lambda(t|x)$. By definition, it also follows that for $t\leq\tau_c(x)$,
\begin{eqnarray*}
\Lambda_n(t|x)-\Lambda(t|x)&=&\int_{-\infty}^t\dfrac{dH_n^u(s|x)}{1-H_n(s^-|x)}-\int_{-\infty}^t\dfrac{dH^u(s|x)}{1-H(s^-|x)}\\
&=&\int_{-\infty}^t \Big[\dfrac{1}{1-H_n(s^-|x)}-\dfrac{1}{1-H(s^-|x)} \Big] dH^u(s|x) \\
&& +\int_{-\infty}^t\dfrac{1}{1-H(s^-|x)}d(H_n^u-H^u)(s|x)\\
&&+\int_{-\infty}^t \Big[\dfrac{1}{1-H_n(s^-|x)}-\dfrac{1}{1-H(s^-|x)} \Big] d(H_n^u-H^u)(s|x).
\end{eqnarray*}
Writing the integrand in the first term as
\begin{eqnarray*}
\dfrac{H_n(s^-|x)-H(s^-|x)}{(1-H_n(s^-|x))(1-H(s^-|x))}=\dfrac{H_n(s^-|x)-H(s^-|x)}{(1-H(s^-|x))^2}+\dfrac{(H_n(s^-|x)-H(s^-|x))^2}{(1-H(s^-|x))^2(1-H_n(s^-|x))},
\end{eqnarray*}
and integrating by parts in the second term (see Theorem A.1.2 in \citealt{Fleming1991}), we have by continuity of $H(\cdot|x)$ over $(-\infty,\tau_c(x))$,
\begin{eqnarray*}
&& \int_{-\infty}^t\dfrac{1}{1-H(s^-|x)}d(H_n^u-H^u)(s|x) \\
&& = \dfrac{H_n^u(t|x)-H^u(t|x)}{1-H(t^-|x)} - \int_{-\infty}^t [H_n^u(s|x)-H^u(s|x)] d\left(\dfrac{1}{1-H}\right)(s^-|x)\\
&& = \dfrac{H_n^u(t|x)-H^u(t|x)}{1-H(t^-|x)} - \int_{-\infty}^t\dfrac{H_n^u(s|x)-H^u(s|x)}{(1-H(s^-|x))^2}dH(s^-|x).
\end{eqnarray*}
It now follows that
\begin{eqnarray*}
\Lambda_n(t|x)-\Lambda(t|x)&=&\int_{-\infty}^t\dfrac{H_n(s^-|x)-H(s^-|x)}{(1-H(s^-|x))^2}dH^u(s|x)+\dfrac{H_n^u(t|x)-H^u(t|x)}{1-H(t^-|x)}\\
\nonumber&&-\int_{-\infty}^t\dfrac{H_n^u(s|x)-H^u(s|x)}{(1-H(s^-|x))^2}dH(s^-|x)+R_{n,2}(t|x)+R_{n,3}(t|x),
\end{eqnarray*}
where
\begin{eqnarray*}
R_{n,2}(t|x)&=&\int_{-\infty}^t\dfrac{(H_n(s^-|x)-H(s^-|x))^2}{(1-H(s^-|x))^2(1-H_n(s^-|x))}dH^u(s|x),\\
R_{n,3}(t|x)&=&\int_{-\infty}^t \Big[\dfrac{1}{1-H_n(s^-|x)}-\dfrac{1}{1-H(s^-|x)} \Big] d(H_n^u-H^u)(s|x).
\end{eqnarray*}
By Lemma \ref{lemma_rates} we have $H_n(\tau_c(x)^-|x)\overset{\mathbb{P}}{\rightarrow} H(\tau_c(x)^-|x)<1$ as $n\rightarrow \infty$, meaning that we may suppose that for $n$ large enough, $H_n(\tau_c(x)^-|x)<1$ (in $\mathbb{P}$). This gives
\begin{eqnarray*}
\sup_{t\leq\tau_c(x)}|R_{n,2}(t|x)|&\leq&\sup_{t\in\mathbb{R}}|H_n(t|x)-H(t|x)|^2\dfrac{1}{(1-H_n(\tau_c(x)^-|x)(1-H(\tau_c(x)^-|x))^2}\\
&=&\mathcal{O}_\mathbb{P}((nh^p)^{-1} |\log h|).
\end{eqnarray*}
In order to uniformly bound the term $R_{n,3}(t|x)$, we use the proof of Lemma 2.1 in \cite{VanKeilegom1997}. Hence, we define $\tau<\tau_c(x)$ such that $H(\tau|x)<\varepsilon$ for $\varepsilon>0$ along with a partition of the interval $[\tau,\tau_c(x)]$ into $k_n=O((nh^p)^{1/2} |\log h|^{-1/2})$ subintervals $(t_i,t_{i+1}]$ of length $k_n=O((nh^p)^{-1/2} |\log h|^{1/2})$. We have that 
\begin{eqnarray}
\sup_{\tau\leq t\leq\tau_c(x)}|R_{n,3}(t|x)|&\leq& 2 \max_{1\leq i\leq k_n}\sup_{s\in(t_i,t_{i+1}]}\left\vert\dfrac{1}{H_n(s^-|x)}-\dfrac{1}{H_n(t_i^-|x)}-\dfrac{1}{H(s^-|x)}+\dfrac{1}{H(t_i^-|x)}\right\vert \nonumber \\[.2cm]
&&+\,k_n\dfrac{\sup_{0\leq t\leq \tau_c(x)}\vert H_n(t|x)-H(t|x)\vert}{(1-H(\tau_c(x)^-|x))(1-H_n(\tau_c(x)^-|x))} \nonumber \\[.2cm]
&&\,\times \max_{1\leq i\leq k_n}\left\vert H_n^u(t_{i+1}|x)-H_n^u(t_i|x) - H^u(t_{i+1}|x) + H^u(t_i|x)\right\vert. \label{modcont}
\end{eqnarray}
Since the functions $H(\cdot^-|x)$ and $H^u(\cdot|x)$ are continuous, the same arguments as in the proof of Lemma 2.1 in \cite{VanKeilegom1997} show that
\begin{eqnarray*}
\sup_{\tau\leq t\leq\tau_c(x)}|R_{n,3}(t|x)|=\mathcal{O}_\mathbb{P}((nh^p)^{-3/4} |\log h|^{3/4})
\end{eqnarray*}
for any $\tau$.
From those bounds, one can deduce that
\begin{eqnarray*}
\sup_{\tau\leq t\leq\tau_c(x)}|\Lambda_n(t|x)-\Lambda(t|x)|=\mathcal{O}_\mathbb{P}((nh^p)^{-1/2}|\log h|^{1/2}),
\end{eqnarray*}
giving that
\begin{eqnarray*}
\sup_{\tau\leq t\leq\tau_c(x)}|R_{n,1}(t|x)|=\mathcal{O}_\mathbb{P}((nh^p)^{-1}|\log h|),
\end{eqnarray*}
which shows our result.\hfill$\Box$

\bigskip

\noindent
{\bf Proof of Theorem \ref{theorem_beran}.} According to Theorem \ref{proposition_iid}, the asymptotic properties of the empirical process $(nh^p)^{1/2}(F_n(\cdot|x)-F(\cdot|x))$ are determined by those of the process
\begin{eqnarray}
\label{proof_process}
\left\{(nh^p)^{1/2}\sum_{i=1}^nW_h(x-X_i)g(t,T_i,\delta_i|x),t\in[\tau,\tau_c(x)]\right\},
\end{eqnarray}
since $\sup_{\tau\leq t\leq\tau_c(x)}(nh^p)^{1/2}|r_n(t|x)|=\mathcal{O}_\mathbb{P}\left((nh^p)^{-1/4} |\log h|^{3/4}\right) = o_\mathbb{P}(1)$. Next, note that
\begin{eqnarray*}
&&\sum_{i=1}^nW_h(x-X_i)g(t,T_i,\delta_i|x)\\
&&=\dfrac{1}{f(x)}\dfrac{1}{n}\sum_{i=1}^nK_h(x-X_i) g(t,T_i,\delta_i|x)\\
&&\quad-\quad\dfrac{1}{n}\sum_{i=1}^nK_h(x-X_i)g(t,T_i,\delta_i|x)\dfrac{1}{f(x)^2}\left(\dfrac{1}{n}\sum_{i=1}^nK_h(x-X_i)-f(x)\right)(1+o_\mathbb{P}(1))\\
&&=\dfrac{1}{f(x)}\left(\dfrac{1}{n}\sum_{i=1}^nK_h(x-X_i)g(t,T_i,\delta_i|x)-\mathbb{E}\left[K_h(x-X)g(t,T,\delta|x)\right]\right)\\
&&\quad\times\left(1-\dfrac{1}{f(x)}\left(\dfrac{1}{n}\sum_{i=1}^nK_h(x-X_i)-f(x)\right)(1+o_\mathbb{P}(1))\right)\\
&&\quad+\dfrac{1}{f(x)}\mathbb{E}\left[K_h(x-X)g(t,T,\delta|x)\right]\left(1-\dfrac{1}{f(x)}\left(\dfrac{1}{n}\sum_{i=1}^nK_h(x-X_i)-f(x)\right)(1+o_\mathbb{P}(1))\right),
\end{eqnarray*}
and hence, the weak convergence of the process (\ref{proof_process}) is equivalent to that of the process
\begin{eqnarray}
\label{proof_process2}
&& \left\{(nh^p)^{1/2} \dfrac{1}{f(x)} \left(\dfrac{1}{n}\sum_{i=1}^nK_h(x-X_i)g(t,T_i,\delta_i|x)-\mathbb{E}\left[K_h(x-X)g(t,T,\delta|x)\right]\right) \right., \nonumber \\
&& \hspace*{10cm} \frac{}{} t\in[\tau,\tau_c(x)]\Big\},
\end{eqnarray}
since by (\ref{proof_rate1}) and (\ref{kerneldensity}),

\begin{eqnarray*}
\sup_{\tau\leq t\leq\tau_c(x)} \left|(nh^p)^{1/2}\mathbb{E}\left[K_h(x-X)g(t,T,\delta|x)\right]\right|&=& \mathcal{O}((nh^p)^{1/2} h^\eta) = o(1),\\
\dfrac{1}{n}\sum_{i=1}^nK_h(x-X_i)-f(x) &=& o_\mathbb{P}(1).
\end{eqnarray*}

Using integration by parts, one can rewrite the function $g(t,T,\delta|x)$ as
\begin{eqnarray*}
\label{proof_g}
g(t,T,\delta|x)=(1-F(t|x))\left\{\dfrac{\ind_{\{\delta=1,T\leq t\}}}{1-H(T^-|x)}-\int_{-\infty}^{T\wedge t}\dfrac{dH^u(s|x)}{(1-H(s^-|x))^2}\right\}.
\end{eqnarray*}

In order to show the convergence of the stochastic process  (\ref{proof_process2}), we will use Theorem 19.28 in \cite{Vaart1998}. To do so, we need to introduce some more notations. Let $P$ denote the law of the vector $(T,\delta,X)$ and define the expectation under $P$, its empirical version and the empirical process as follows :
\begin{eqnarray*}
Pf=\int fdP,\hspace{0.5cm}\mathbb{P}_nf=\dfrac{1}{n}\sum_{i=1}^nf\left(T_i, \delta_i, X_i\right),\hspace{0.5cm}\mathbb{G}_nf=\sqrt{n}(\mathbb{P}_n-P)f,
\end{eqnarray*}
for any real-valued measurable function $f$. We also introduce our sequence of classes $\mathcal{F}_n$ with functions taking values in $E=\mathbb{R}\times\{0,1\}\times\mathbb{R}^p$ as
\begin{eqnarray*}
\mathcal{F}_n &=& \left\{(u,v,w)\rightarrow f_{n,t}(u,v,w), \, t\in[\tau,\tau_c(x)]\right\} \\
&=& \left\{(u,v,w)\rightarrow\sqrt{h^p}K_h(x-w)g(t,u,v|x),\, t\in[\tau,\tau_c(x)]\right\}.
\end{eqnarray*}
Denote now by $E_n$ an envelope function of the class $\mathcal{F}_n$. According to Theorem 19.28 in \cite{Vaart1998}, the weak convergence of the stochastic process (\ref{proof_process2}) follows from the following four conditions :
\begin{eqnarray}
\label{ST1}
\sup_{\rho(t,s)\leq \delta_n} P(f_{n,t}-f_{n,s})^2 &\longrightarrow& 0 \mbox{ for every $\delta_n \searrow 0$,}\\
\label{ST2}
PE_n^2&=&O(1),\\
\label{ST3}
PE_n^2\{E_n>\varepsilon \sqrt n\} &\longrightarrow& 0 \mbox{ for every $\varepsilon >0,$}\\
\label{ST4}
J(\delta_n, {\cal F}_n, L_2)&\longrightarrow& 0 \mbox{ for every $\delta_n \searrow 0$,}
\end{eqnarray}
where $\rho$ is a semimetric that makes $[\tau,\tau_c(x)]$ a totally bounded space.   We will work with $\rho(t) = |t|$ for any $t \in \mathbb{R}$.  For an appropriate constant $M>0$,
\begin{eqnarray*}
E_n(u,v,w)=\sqrt{h^p}K_h(x-w)M,
\end{eqnarray*}
since $g(t,T,\delta|x)$ is uniformly bounded for $t\in[\tau,\tau_c(x)]$.

We start by proving $(\ref{ST1})$. By definition we have for $n$ large enough,
\begin{eqnarray*}
P(f_{n,t}-f_{n,s})^2\leq\int_{S_K}K^2(u)\mathbb{E}[(g(t,T,\delta|x)-g(s,T,\delta|x))^2|X=x-hu]f(x-hu)du.
\end{eqnarray*}
In particular, for $s<t\leq\tau_c(x)$, 
\begin{eqnarray*}
&& |g(t,T,\delta|x)-g(s,T,\delta|x)| \\
&& \leq |F(t|x)-F(s|x)|\left|\dfrac{\ind_{\{\delta=1,T\leq t\}}}{1-H(T^-|x)}-\int_{-\infty}^{T\wedge t}\dfrac{dH^u(s|x)}{(1-H(s^-|x))^2}\right|\\
&& \quad + (1-F(s|x))\left|\dfrac{\ind_{\{\delta=1,s<T\leq t\}}}{1-H(T^-|x)}-\int_{T\wedge s}^{T\wedge t}\dfrac{dH^u(s|x)}{(1-H(s^-|x))^2}\right|\\
&& \le |F(t|x)-F(s|x)|\left[\dfrac{1}{1-H(\tau_c(x)^-|x)}+\int_{-\infty}^{\tau_c(x)}\dfrac{dH^u(s|x)}{(1-H(s^-|x))^2}\right]\\
&& \quad + \left|\dfrac{\ind_{\{\delta=1,s<T\leq t\}}}{1-H(\tau_c(x)^-|x)}\right|+\left|\int_{T\wedge s}^{T\wedge t}\dfrac{dH^u(s|x)}{(1-H(s^-|x))^2}\right|.
\end{eqnarray*}
Hence, one can find a positive constant $M'>0$ such that
\begin{eqnarray*}
&& |g(t,T,\delta|x)-g(s,T,\delta|x)| \\
&& \leq M'\left[|F(t|x)-F(s|x)|+\ind_{\{\delta_i=1,s<T_i\leq t\}}+\left|H^u(T\wedge t|x)-H^u(T\wedge s|x)\right|\right].
\end{eqnarray*}
Furthermore, we have 
\begin{eqnarray*}
&& \left|H^u(T\wedge t|x)-H^u(T\wedge s|x)\right| \\
&& = \left|H^u(T|x)-H^u(s|x)\right|\ind_{\{s<T\leq t\}} + \left|H^u(t|x)-H^u(s|x)\right|\ind_{\{t<T\}}\\
&& \leq 2\left|H^u(t|x)-H^u(s|x)\right|.
\end{eqnarray*}
This gives
\begin{eqnarray*}
&&\mathbb{E}[(g(t,T,\delta|x)-g(s,T,\delta|x))^2|X=x-hu]f(x-hu)\\
&& \leq (M')^2\left[|F(t|x)-F(s|x)|+\mathbb{E}[\ind_{\{\delta_i=1,s<T_i\leq t\}}|X=x-hu] \right.\\
&& \quad \left. +2\left|H^u(t|x)-H^u(s|x)\right|\right]^2f(x-hu)\\
&&\leq 2(M')^2f(x)\left[|F(t|x)-F(s|x)|^2+2|H^u(t|x-hu)-H^u(s|x-hu)| \right.\\
&& \quad \left. +4\left|H^u(t|x)-H^u(s|x)\right|^2\right] +14 (M')^2\left|f(x)-f(x-hu)\right|,
\end{eqnarray*}
and thus we obtain
\begin{eqnarray*}
P(f_{n,t}-f_{n,s})^2&\leq& 2(M')^2\Vert K\Vert^2_2f(x)\left[|F(t|x)-F(s|x)|^2+2\left|H^u(t|x)-H^u(s|x)\right|\right.\\
&&+\left.4\left|H^u(t|x)-H^u(s|x)\right|^2\right]+ h^\eta C \int_{S_K}\Vert u\Vert^\eta K(u)^2du,
\end{eqnarray*}
for some $C<\infty$.  It follows that $P(f_{n,t}-f_{n,s})^2$ converges uniformly in $(t,s)$ towards 0 as $|t-s|\to 0$ and $n \to \infty$, by the uniform continuity of the functions $F(\cdot|x)$ and $H^u(\cdot|x)$ over the compact set $[\tau,\tau_c(x)]$.

Now, we move to the proof of (\ref{ST2}) and (\ref{ST3}). It follows that
\begin{eqnarray*}
&&PE_n^2=M^2\int_{S_K}K^2(u)f(x-hu)du\leq h^\eta M^2\int_{S_k}\Vert u\Vert^\eta K(u)^2du+M^2f(x)\Vert K\Vert^2_2,\\
&&PE_n^2\{E_n>\varepsilon\sqrt{n}\}\leq M^2\int_{\{K(u)>M^{-1}\varepsilon\sqrt{nh^p}\}}K^2(u)f(x-hu)du=0,
\end{eqnarray*}
for all $\varepsilon>0$ and $n$ sufficiently large, since $nh^p\rightarrow \infty$ and $K$ is bounded.

Finally, it remains to prove (\ref{ST4}). It is clear that the function classes $\{(u,v,w)\to t,\,t\in[\tau,\tau_c(x)]\}$, $\{(u,v,w)\to u\}$ and $\{(u,v,w)\to\ind_{\{v=1,u\leq t\}},\,t\in[\tau,\tau_c(x)]\}$ are \textit{VC}. Invoking Lemma 2.6.18 {\it(i)}, {\it(vi)} and {\it(viii)} in \cite{Vaart1996}, we obtain that the following classes are also \textit{VC} :
\begin{eqnarray*}
\mathcal{G}_1&=&\left\{(u,v,w)\to\dfrac{\ind_{\{v=1,u\leq t\}}}{1-H(u^-|x)},\,t\in[\tau,\tau_c(x)]\right\},\\
\mathcal{G}_2&=&\left\{(u,v,w)\to\int_{-\infty}^{u\wedge t}\dfrac{dH^u(s|x)}{(1-H(s^-|x))^2},\,t\in[\tau,\tau_c(x)]\right\}.
\end{eqnarray*}
Following the same argument as in the proof of Lemma \ref{lemma_rates}, we can add the kernel function $K$ to the previous function classes and thus obtain that
\begin{eqnarray*}
\mathcal{G}_{n,1}&=&\left\{(u,v,w)\to\sqrt{h^p}K_h(x-w)\dfrac{\ind_{\{v=1,u\leq t\}}}{1-H(u^-|x)},\,t\in[\tau,\tau_c(x)]\right\},\\
\mathcal{G}_{n,2}&=&\left\{(u,v,w)\to\sqrt{h^p}K_h(x-w)\int_{-\infty}^{u\wedge t}\dfrac{dH^u(s|x)}{(1-H(s^-|x))^2},\,t\in[\tau,\tau_c(x)]\right\},
\end{eqnarray*} 
are respectively \textit{VC} with constants $A_i$ and $\nu_i$, $i=1,2$, not depending on $n$ and with shared envelope function $E_n$. Finally, since our class of interest $\mathcal{F}_n$ is included in the class of functions $
\widetilde {\mathcal{F}}_n=\mathcal{G}_{n,1}+\mathcal{G}_{n,2}$ with envelope function $2E_n$, using Lemma 16 in \cite{Nolan1987}, we have for any $t>0$,
\begin{eqnarray*}
\sup_\mathcal{Q} N(\mathcal{F}_n,L_2(Q),2t\Vert E_n\Vert_{Q,2})&\leq&\sup_\mathcal{Q} N(\widetilde{\mathcal{F}}_n,L_2(Q),2t\Vert E_n\Vert_{Q,2})\\
&\leq& \left(\dfrac{2A_1}{t}\right)^{\nu_1}\left(\dfrac{2A_2}{t}\right)^{\nu_2} \leq L\left(\dfrac{1}{\tau}\right)^V,
\end{eqnarray*}
for some $L$ and $V$. Thus, (\ref{ST4}) is established since for any sequence $\delta_n\searrow 0$ and $n$ large enough, we have
\begin{eqnarray*}
J(\delta_n,\mathcal{F}_n,L_2)
&\leq&\int_0^{\delta_n}\sqrt{\log(2^VL)-V\log(t)}dt=o(1).
\end{eqnarray*}
This achieves the proof of the weak convergence, since the covariance structure follows from the proof of Lemma A2 in \cite{VanKeilegom1997}. Finally, we establish the continuity of the process thanks to a sufficient condition due to \cite{Fernique1964}. Indeed, let $(s,t)\in\mathbb{R}^2$ and denote $\widebar{F}(\cdot|x)=1-F(\cdot|x)$.  Then,
\begin{eqnarray*}
&&\mathbb{E}[(Z(s|x)-Z(t|x))^2]\dfrac{f(x)}{\Vert K\Vert ^2_2}\\
&&=\widebar{F}(s|x)^2\int_{-\infty}^{s}\dfrac{dH^u(y|x)}{(1-H(y^-|x))^2} -2\widebar{F}(s|x)\widebar{F}(t|x)\int_{-\infty}^{s\wedge t}\dfrac{dH^u(y|x)}{(1-H(y^-|x))^2} \\
&& \quad + \widebar{F}(t|x)^2\int_{-\infty}^{t}\dfrac{dH^u(y|x)}{(1-H(y^-|x))^2}\\
&&=\widebar{F}(s|x)\left[\widebar{F}(s|x)\int_{s\wedge t}^{s}\dfrac{dH^u(y|x)}{(1-H(y^-|x))^2}+\int_{-\infty}^{s\wedge t}\dfrac{dH^u(y|x)}{(1-H(y^-|x))^2}(F(t|x)-F(s|x))\right]\\
&& \quad +\widebar{F}(t|x)\left[\widebar{F}(t|x)\int_{s\wedge t}^{t}\dfrac{dH^u(y|x)}{(1-H(y^-|x))^2}+\int_{-\infty}^{s\wedge t}\dfrac{dH^u(y|x)}{(1-H(y^-|x))^2}(F(s|x)-F(t|x))\right]\\
&&=\widebar{F}(s\vee t|x)^2\int_{s\wedge t}^{s\vee t}\dfrac{dH^u(y|x)}{(1-H(y^-|x))^2}+\int_{-\infty}^{s\wedge t}\dfrac{dH^u(y|x)}{(1-H(y^-|x))^2}(F(s|x)-F(t|x))^2\\
&&\leq \dfrac{c|s-t|^{\eta'}}{(1-H(\tau_c(x)^-|x))^2}.
\end{eqnarray*}
This yields that $\sqrt{\mathbb{E}[(Z(s|x)-Z(t|x))^2]}\leq\xi(t-s)$, with 
$$ \xi(t-s)=\sqrt{\dfrac{c}{f(x)}}\dfrac{\Vert K\Vert_2|s-t|^{\eta'/2}}{1-H(\tau_c(x)^-|x)} $$ 
being monotone and
\begin{eqnarray*}
\int_0^1\dfrac{\xi(u)}{u|\log(u)|^{1/2}}du<+\infty. \\[-1.6cm]
\end{eqnarray*}
\hfill$\Box$\\

\subsection{Proofs of Section \ref{ours}}

\begin{lemma}  \label{lem1}
Assume that $\tau_c=\tau_c(x')$ does not depend on $x'\in S_X$. Then, under the conditions of Theorem \ref{proposition_gamma},
\begin{eqnarray}
\label{scaledprocess}
F_n(t\tau_n|x)-F(t\tau_c|x) = Z_n(t\tau_c|x) + r_n(t\tau_c|x),
\end{eqnarray}
where $\sup_{0 \le t \le 1} |r_n(t\tau_c|x)| = o_\mathbb{P}((nh^p)^{-1/2})$ and 
$$ Z_n(t\tau_c|x) = \sum_{i=1}^n W_h(x-X_i) g(t\tau_c,T_i,\delta_i|x), $$
where $g(t\tau_c,T_i,\delta_i|x)$ is defined in the statement of Theorem \ref{proposition_iid}.
\end{lemma}

\noindent
{\bf Proof.} 
Write
$$ F_n(t\tau_n|x)-F(t\tau_c|x) = [F_n(t\tau_n|x)-F(t\tau_n|x)] + [F(t\tau_n|x)-F(t\tau_c|x)] = A(t|x) + B(t|x). $$
For the term $B(t|x)$, note that by Assumption $(\mathcal{H}.4)$, for $C'>0$ large enough,
\begin{eqnarray*}
\left|F(t\tau_n|x)-F(t\tau_c|x)\right|\leq C'\left|\tau_n-\tau_c\right|^{\eta'}.
\end{eqnarray*} 
We will prove that $(nh^p)^{1/2}\left|\tau_n-\tau_c\right|^{\eta'}=o_\mathbb{P}(1)$. Note that for $\varepsilon>0$ we have that
\begin{eqnarray*}
\mathbb{P}((nh^p)^{1/2}\left|\tau_n-\tau_c\right|^{\eta'}\geq\varepsilon)&\leq&\mathbb{P}(\tau_n<\tau_c)=[\mathbb{P}(T<\tau_c)]^n,
\end{eqnarray*}
and by definition,
\begin{eqnarray*}
\mathbb{P}(T<\tau_c)&=&\int_{S_X}f(u)\mathbb{P}(T<\tau_c|X=u)du\\
&=&1-\int_{S_X}f(u)\mathbb{P}(C\geq\tau_c|X=u)\mathbb{P}(Y\geq\tau_c|X=u)du\\
&\leq&1-\inf_{x'\in B(x,\delta)} [G(\tau_c|x')-G(\tau_c^-|x')] \int_{B(x,\delta)}f(u)\mathbb{P}(Y\geq\tau_c|X=u)du<1
\end{eqnarray*}
by assumption $(\mathcal{A}')$, implying that $[\mathbb{P}(T<\tau_c)]^n\to 0$ as $n\to\infty$.

Next, it follows from Theorem \ref{proposition_iid} that $A(t|x) = Z_n(t\tau_n|x) + o_\mathbb{P}((nh^p)^{-1/2})$ uniformly in $0 \le t \le 1$.   Note that we can write $Z_n(t\tau_n|x)$ as
\begin{align*}
& Z_n(t\tau_n|x) \\
& = (1-F(t \tau_n|x))  \left\{\int_{-\infty}^{t\tau_n} \dfrac{H_n(s^-|x)-H(s^-|x)}{(1-H(s^-|x))^2}dH^u(s|x)+\dfrac{H_n^u(t\tau_n|x)-H^u(t\tau_n|x)}{1-H(t\tau_n^-|x)} \right. \\
& \quad \left. -\int_{-\infty}^{t\tau_n}\dfrac{H_n^u(s|x)-H^u(s|x)}{(1-H(s^-|x))^2}dH(s^-|x)\right\} \\
& = (1-F(t \tau_c|x))  \left\{\int_{-\infty}^{t\tau_c} \dfrac{H_n(s^-|x)-H(s^-|x)}{(1-H(s^-|x))^2}dH^u(s|x)+\dfrac{H_n^u(t\tau_c|x)-H^u(t\tau_c|x)}{1-H(t\tau_c^-|x)} \right. \\
& \quad \left. -\int_{-\infty}^{t\tau_c}\dfrac{H_n^u(s|x)-H^u(s|x)}{(1-H(s^-|x))^2}dH(s^-|x)\right\} + o_\mathbb{P}((nh^p)^{-1/2}) \\
& = Z_n(t\tau_c|x) + o_\mathbb{P}((nh^p)^{-1/2}),
\end{align*}
where the second equality above follows from Assumptions $(\mathcal{A})$ and $(\mathcal{H})$, Lemma \ref{lemma_rates}, the modulus of continuity of the estimator $H_n^u$ given in (\ref{modcont}), and from the fact that $\tau_n-\tau_c = o_\mathbb{P}((nh^p)^{-1/(2\eta')})$.\hfill $\Box$

\bigskip

\noindent
\textbf{Proof of Theorem \ref{proposition_gamma}.}  The proof is mainly based on the  representation in (\ref{scaledprocess}). Hence, for $n$ large enough and thanks to straightforward Taylor expansions, we have
\begin{eqnarray*}
&&\dfrac{F_n(y_2^ 2\tau_n|x)-F_n(y_2\tau_n|x)}{F_n(y_2\tau_n|x)-F_n(\tau_n|x)}
\\
&&=\dfrac{F(y_2^ 2\tau_c|x)-F(y_2\tau_c|x)}{F(y_2\tau_c|x)-F(\tau_c|x)}
+\dfrac{F_n(y_2^ 2\tau_n|x)-F(y_2^2\tau_c|x)-F_n(y_2\tau_n|x)+F(y_2\tau_c|x)}{F(y_2\tau_c|x)-F(\tau_c|x)}\\
&& \quad + (F_n(y_2^ 2\tau_n|x)-F_n(y_2\tau_n|x))\left[\dfrac{1}{F_n(y_2\tau_n|x)-F_n(\tau_n|x)}-\dfrac{1}{F(y_2\tau_c|x)-F(\tau_c|x)}\right]\\
&&=\dfrac{F(y_2^ 2\tau_c|x)-F(y_2\tau_c|x)}{F(y_2\tau_c|x)-F(\tau_c|x)}+\dfrac{Z_n(y_2^2\tau_c|x)-Z_n(y_2\tau_c|x)}{F(y_2\tau_c|x)-F(\tau_c|x)}+o_\mathbb{P}((nh^p)^{-1/2})\\
&& \quad+(F(y_2^ 2\tau_c|x)-F(y_2\tau_c|x)+o_\mathbb{P}(1))\dfrac{Z_n(\tau_c|x)-Z_n(y_2\tau_c|x)}{(F(y_2\tau_c|x)-F(\tau_c|x))^2}(1+o_\mathbb{P}(1))\\
&&=\dfrac{F(y_2^ 2\tau_c|x)-F(y_2\tau_c|x)}{F(y_2\tau_c|x)-F(\tau_c|x)}+\dfrac{Z_n(y_2^2\tau_c|x)-Z_n(y_2\tau_c|x)}{F(y_2\tau_c|x)-F(\tau_c|x)}\\
&& \quad +\dfrac{F(y_2^ 2\tau_c|x)-F(y_2\tau_c|x)}{(F(y_2\tau_c|x)-F(\tau_c|x))^2}(Z_n(\tau_c|x)-Z_n(y_2\tau_c|x))+o_\mathbb{P}((nh^p)^{-1/2}).
\end{eqnarray*}
Next, we use the Delta method applied to the function $x\rightarrow-1/\log_{y_2}(x)$ and obtain
\begin{eqnarray}
\label{proof_gamma}
\nonumber
\widehat\gamma(x) \! &=& \! \gamma_{y_2,\tau_c}(x)+\phi_{y_2}\left(\dfrac{F(y_2^ 2\tau_c|x)-F(y_2\tau_c|x)}{F(y_2\tau_c|x)-F(\tau_c|x)}\right)\left\{\dfrac{Z_n(y_2^2\tau_c|x)-Z_n(y_2\tau_c|x)}{F(y_2\tau_c|x)-F(\tau_c|x)}\right.\\
&& +\left.\dfrac{F(y_2^ 2\tau_c|x)-F(y_2\tau_c|x)}{(F(y_2\tau_c|x)-F(\tau_c|x))^2}(Z_n(\tau_c|x)-Z_n(y_2\tau_c|x))+o_\mathbb{P}((nh^p)^{-1/2})\right\}\!.
\end{eqnarray} 
Finally, the asymptotic normality of $(nh^p)^{1/2}(\widehat\gamma(x) - \gamma_{y_2,\tau_c}(x))$ follows by Lyapunov's central limit theorem.  \hfill$\Box$

\bigskip

\noindent\textbf{Proof of Theorem \ref{theorem_p}.} Similarly to the proof of Theorem \ref{proposition_gamma}, we will use the representation in (\ref{scaledprocess}). Write
\begin{align*}
& \widehat p(x) - F_n(\tau_n|x) \\
& =(F_n(\tau_n|x)-F_n(y_1\tau_n|x))\left\{{1\over y_1^{-1/\gamma_{y_2,\tau_c}(x)}-1} +{1\over y_1^{-1/\widehat\gamma(x)}-1}-{1\over y_1^{-1/\gamma_{y_2,\tau_c}(x)}-1}\right\}\\
&= (F_n(\tau_n|x)-F_n(y_1\tau_n|x))\left\{{1\over y_1^{-1/\gamma_{y_2,\tau_c}(x)}-1} +(\widehat\gamma(x)-\gamma_{y_2,\tau_c}(x))\psi_{y_1}(\gamma_{y_2,\tau_c}(x))(1+o_\mathbb{P}(1))\right\},
\end{align*}
where 
\begin{eqnarray*}
\psi_{y_1}(x)=\dfrac{-\log(y_1)y_1^{-1/x}}{(x(y_1^{-1/x}-1))^2}.
\end{eqnarray*}
By (\ref{proof_gamma}), this leads to
\begin{eqnarray*}
\widehat p(x)&=& F(\tau_c|x)+Z_n(\tau_c|x)+o_\mathbb{P}((nh^p)^{-1/2}) \\
&&+ \big(F(\tau_c|x)-F(y_1\tau_c|x)+Z_n(\tau_c|x)-Z_n(y_1\tau_c|x)+o_\mathbb{P}((nh^p)^{-1/2})\big) \\
&&\times \left[{1\over y_1^{-1/\gamma_{y_2,\tau_c}(x)}-1}+\psi_{y_1}(\gamma_{y_2,\tau_c}(x))\phi_{y_2}\left(\dfrac{F(y_2^ 2\tau_c|x)-F(y_2\tau_c|x)}{F(y_2\tau_c|x)-F(\tau_c|x)}\right) \right. \\
&& \quad \left. \times \left\{\dfrac{Z_n(y_2^2\tau_c|x)-Z_n(y_2\tau_c|x)}{F(y_2\tau_c|x)-F(\tau_c|x)} + \dfrac{F(y_2^ 2\tau_c|x)-F(y_2\tau_c|x)}{(F(y_2\tau_c|x)-F(\tau_c|x))^2}(Z_n(\tau_c|x)-Z_n(y_2\tau_c|x)) \right. \right. \\
&& \quad  +o_\mathbb{P}((nh^p)^{-1/2})\Big\}  (1+o_\mathbb{P}(1))\Big].
\end{eqnarray*}
By rearranging the terms, we obtain
\begin{eqnarray}
\label{proof_p_n}
\widehat p(x)&=&p_{y_1,y_2,\tau_c}(x)+Z_n(\tau_c|x)+\dfrac{Z_n(\tau_c|x)-Z_n(y_1\tau_c|x)}{y_1^{-1/\gamma_{y_2,\tau_c}(x)}-1}\\
\nonumber && - \dfrac{F(y_1\tau_c|x)-F(\tau_c|x)}{F(y_2\tau_c|x)-F(\tau_c|x)}\psi_{y_1}(\gamma_{y_2,\tau_c}(x))\phi_{y_2}\left(\dfrac{F(y_2^ 2\tau_c|x)-F(y_2\tau_c|x)}{F(y_2\tau_c|x)-F(\tau_c|x)}\right) \\ 
\nonumber && \times \left[Z_n(y_2^2\tau_c|x)-Z_n(y_2\tau_c|x)+\dfrac{F(y_2^ 2\tau_c|x)-F(y_2\tau_c|x)}{F(y_2\tau_c|x)-F(\tau_c|x)}(Z_n(\tau_c|x)-Z_n(y_2\tau_c|x))\right] \\
\nonumber && +o_\mathbb{P}((nh^p)^{-1/2}).
\end{eqnarray}
The result now follows from Lyapunov's central limit theorem. \hfill$\Box$

\bigskip

\noindent\textbf{Proof of Theorem \ref{theorem_F}.}  We will look separately at the different terms in the formula of $\widehat F(t|x)$.  First note that similarly as in the proof of Lemma \ref{lem1}, we can show that 
\begin{eqnarray}
\label{shiftedprocess}
F_n(t\wedge\tau_n|x)=F(t\wedge\tau_c|x)+Z_n(t\wedge\tau_c|x)+o_\mathbb{P}((nh^p)^{-1/2}),
\end{eqnarray}
uniformly in $t$. 

Next, since the asymptotic representations of the estimators $F_n(\tau_n|x)$ and $\widehat p(x)$ have already been developed (see Lemma \ref{lem1} and Theorem \ref{theorem_p}), we can focus on developing an asymptotic expansion for the process
\begin{eqnarray*}
\left\{\left[\frac{t}{\tau_n}\vee 1\right]^{-1/\widehat\gamma(x)},\,t\in[a,+\infty)\right\}.
\end{eqnarray*}

For any $t\geq a$,
\begin{eqnarray*}
&& \left[\frac{t}{\tau_n}\vee 1\right]^{-1/\widehat\gamma(x)}-\left[\frac{t}{\tau_n}\vee 1\right]^{-1/\gamma_{y_2,\tau_c}(x)} \\
&&=\left(\dfrac{1}{\gamma_{y_2,\tau_c}(x)}-\dfrac{1}{\widehat\gamma(x)}\right)\log\left(\frac{t}{\tau_n}\vee 1\right)\left[\frac{t}{\tau_n}\vee 1\right]^{-1/\gamma_{y_2,\tau_c}(x)}\\
&& \quad +\frac{1}{2}\int_{-1/\gamma_{y_2,\tau_c}(x)}^{-1/\widehat\gamma(x)}\log\left(\frac{t}{\tau_n}\vee 1\right)^2\left[\frac{t}{\tau_n}\vee 1\right]^s \left(s-\dfrac{1}{\gamma_{y_2,\tau_c}(x)}\right)ds.
\end{eqnarray*}
\noindent
Furthermore, one can show that for any $s<0$ and $x\geq 1$, $0\leq\log(x)^2x^s\leq 4e^{-2}s^{-2}$, giving
\begin{eqnarray*}
&& \left|\int_{-1/\gamma_{y_2,\tau_c}(x)}^{-1/\gamma_n}\log\left(\frac{t}{\tau_n}\vee 1\right)^2\left[\frac{t}{\tau_n}\vee 1\right]^s \left(s-\dfrac{1}{\gamma_{y_2,\tau_c}(x)}\right)ds\right| \\
&&\leq 4e^{-2}(\widehat\gamma(x)\wedge\gamma_{y_2,\tau_c}(x))^2 \left(\dfrac{1}{\widehat\gamma(x)}-\dfrac{1}{\gamma_{y_2,\tau_c}(x)}\right)^2 \\
&&=o_\mathbb{P}\left(\left|\dfrac{1}{\widehat\gamma(x)}-\dfrac{1}{\gamma_{y_2,\tau_c}(x)}\right|\right),
\end{eqnarray*}
uniformly in $t\geq a$. It thus follows that
\begin{eqnarray*}
&&\left[\frac{t}{\tau_n}\vee 1\right]^{-1/\widehat\gamma(x)}-\left[\frac{t}{\tau_n}\vee 1\right]^{-1/\gamma_{y_2,\tau_c}(x)}\\
&& = \dfrac{\widehat\gamma(x)-\gamma_{y_2,\tau_c}(x)}{\gamma_{y_2,\tau_c}(x)^2}\log\left(\frac{t}{\tau_n}\vee 1\right)\left[\frac{t}{\tau_n}\vee 1\right]^{-1/\gamma_{y_2,\tau_c}(x)}(1+o_\mathbb{P}(1)),
\end{eqnarray*}
uniformly in $t\geq a$. 

For the terms depending on $t/\tau_n$, we can replace them by $t/\tau_c$ by making use of the fast rate of convergence of $\tau_n$.   In particular, we have that for some positive constant $M$,
\begin{eqnarray*}
&&\left|\left[\frac{t}{\tau_n}\vee 1\right]^{-1/\gamma_{y_2,\tau_c}(x)}-\left[\frac{t}{\tau_c}\vee 1\right]^{-1/\gamma_{y_2,\tau_c}(x)}\right|\leq M|\tau_n-\tau_c|=o_\mathbb{P}((nh^p)^{-1/2}),\\
&&\left|\log\left(\frac{t}{\tau_n}\vee 1\right)\left[\frac{t}{\tau_n}\vee 1\right]^{-1/\gamma_{y_2,\tau_c}(x)}-\log\left(\frac{t}{\tau_c}\vee 1\right)\left[\frac{t}{\tau_c}\vee 1\right]^{-1/\gamma_{y_2,\tau_c}(x)}\right| \\
&& \quad \leq M|\tau_n-\tau_c|=o_\mathbb{P}((nh^p)^{-1/2}).
\end{eqnarray*}
\noindent
This leads to the following representation :
\begin{eqnarray}
\label{proof_t}
&&\left[\frac{t}{\tau_n}\vee 1\right]^{-1/\widehat\gamma(x)}\\
\nonumber&&=\left[\frac{t}{\tau_c}\vee 1\right]^{-1/\gamma_{y_2,\tau_c}(x)}+\dfrac{\widehat\gamma(x)-\gamma_{y_2,\tau_c}(x)}{\gamma_{y_2,\tau_c}(x)^2}\log\left(\frac{t}{\tau_c}\vee 1\right)\left[\frac{t}{\tau_c}\vee 1\right]^{-1/\gamma_{y_2,\tau_c}(x)}(1+o_\mathbb{P}(1)),
\end{eqnarray}
where $\widehat\gamma(x)-\gamma_{y_2,\tau_c}(x)$ is given in (\ref{proof_gamma}).

Thanks to equations (\ref{proof_gamma}), (\ref{proof_p_n}), (\ref{shiftedprocess}) and (\ref{proof_t}), we now have all the ingredients for the asymptotic representation of the process $\widehat{F}(t|x)$. Indeed, for $t\geq a$, we have 
\begin{eqnarray*}
&& \widehat{F}(t|x) \\
&& = F_n(t\wedge\tau_n|x)+\left(\widehat p(x)-F_n(\tau_n|x)\right)\left(1-\left[\frac{t}{\tau_n}\vee 1\right]^{-1/\widehat\gamma(x)}\right)\\
&&= F(t\wedge\tau_c|x)+F_n(t\wedge\tau_n|x)-F(t\wedge\tau_c|x) \\
&& \quad +\left(p_{y_1,y_2,\tau_c}(x)-F(\tau_c|x)+\widehat p(x)-p_{y_1,y_2,\tau_c}(x)-F_n(\tau_n|x)+F(\tau_c|x)\right)\\
&& \quad \times\left(1-\left[\frac{t}{\tau_c}\vee 1\right]^{-1/\gamma_{y_2,\tau_c}(x)}-\dfrac{\widehat\gamma(x)-\gamma_{y_2,\tau_c}(x)}{\gamma_{y_2,\tau_c}(x)^2}\log\left(\frac{t}{\tau_c}\vee 1\right)\left[\frac{t}{\tau_c}\vee 1\right]^{-1/\gamma_{y_2,\tau_c}(x)}(1+o_\mathbb{P}(1))\right)\\
&&=F_{y_1,y_2,\tau_c}(t|x)+F_n(t\wedge\tau_n|x)-F(t\wedge\tau_c|x)\\
&& \quad - \dfrac{p_{y_1,y_2,\tau_c}(x)-F(\tau_c|x)}{\gamma_{y_2,\tau_c}(x)^2}\log\left(\frac{t}{\tau_c}\vee 1\right)\left[\frac{t}{\tau_c}\vee 1\right]^{-1/\gamma_{y_2,\tau_c}(x)}(\widehat\gamma(x)-\gamma_{y_2,\tau_c}(x))\\
&&\quad + \left(1-\left[\frac{t}{\tau_c}\vee 1\right]^{-1/\gamma_{y_2,\tau_c}(x)}\right)(\widehat p(x)-p_{y_1,y_2,\tau_c}(x)-F_n(\tau_n|x)+F(\tau_c|x))+o_\mathbb{P}((nh^p)^{-1/2}).
\end{eqnarray*}

Recall the notations $a(y_2|x)$, $b(y_2|x)$ and $c(y_1,y_2|x)$ introduced before.  Then, 
\begin{eqnarray*}
&& \widehat{F}(t|x)-F_{y_1,y_2,\tau_c}(t|x) \\
&& = Z_n(t\wedge\tau_c|x)\\
&& \quad + Z_n(\tau_c|x)\left[-\dfrac{p(x)-F(\tau_c|x)}{\gamma^2_{y_2,\tau_c}(x)}\log\left(\frac{t}{\tau_c}\vee 1\right)\left[\frac{t}{\tau_c}\vee 1\right]^{-1/\gamma_{y_2,\tau_c}(x)}\phi_{y_2}\left(b(y_2|x)\right)\dfrac{b(y_2|x)}{a(y_2|x)}\right.\\
&& \quad +\left.\left(1-\left[\frac{t}{\tau_c}\vee 1\right]^{-1/\gamma_{y_2,\tau_c}(x)}\right)\left(\dfrac{1}{y^{-1/\gamma_{y_2,\tau_c}(x)}-1}-c(y_1,y_2|x)\psi_{y_1}(\gamma_{y_2,\tau_c}(x))\phi_{y_2}(b(y_2|x))b(y_2|x)\right)\right]\\
&& \quad + Z_n(y_2\tau_c|x)\left[\dfrac{p(x)-F(\tau_c|x)}{\gamma^2_{y_2,\tau_c}(x)}\log\left(\frac{t}{\tau_c}\vee 1\right)\left[\frac{t}{\tau_c}\vee 1\right]^{-1/\gamma_{y_2,\tau_c}(x)}\phi_{y_2}\left(b(y_2|x)\right)\dfrac{1+b(y_2|x)}{a(y_2|x)}\right.\\
&& \quad +\left.\left(1-\left[\frac{t}{\tau_c}\vee 1\right]^{-1/\gamma_{y_2,\tau_c}(x)}\right)c(y_1,y_2|x)\psi_{y_1}(\gamma_{y_2,\tau_c}(x))\phi_{y_2}(b(y_2|x))(1-b(y_2|x))\right]\\
&& \quad - Z_n(y_2^2\tau_c|x)\left[\dfrac{p(x)-F(\tau_c|x)}{\gamma^2_{y_2,\tau_c}(x)}\log\left(\frac{t}{\tau_c}\vee 1\right)\left[\frac{t}{\tau_c}\vee 1\right]^{-1/\gamma_{y_2,\tau_c}(x)}\phi_{y_2}\left(b(y_2|x)\right)\dfrac{1}{a(y_2|x)}\right.\\
&& \quad +\left.\left(1-\left[\frac{t}{\tau_c}\vee 1\right]^{-1/\gamma_{y_2,\tau_c}(x)}\right)c(y_1,y_2|x)\psi_{y_1}(\gamma_{y_2,\tau_c}(x))\phi_{y_2}(b(y_2|x))\right]\\
&& \quad - Z_n(y_1\tau_c|x)\dfrac{1-\left[\frac{t}{\tau_c}\vee 1\right]^{-1/\gamma_{y_2,\tau_c}(x)}}{y^{-1/\gamma_{y_2}(x)}-1}.
\end{eqnarray*}
Finally, to show the weak convergence of the above process, one can follow similar arguments as in the proof of Theorem \ref{theorem_beran}.  This finishes the proof. \hfill $\Box$

\bibliographystyle{apalike}
\bibliography{bibli.bib}

\begin{thebibliography}{}

\bibitem[Amico and Van~Keilegom, 2018]{Amico2018}
Amico, M. and Van~Keilegom, I. (2018).
\newblock Cure models in survival analysis.
\newblock {\em Annual Review of Statistics and its Application}, 5:311--342.

\bibitem[Beirlant and Guillou, 2001]{Beirlant2001}
Beirlant, J. and Guillou, A. (2001).
\newblock Pareto index estimation under moderate right censoring.
\newblock {\em Scandinavian Actuarial Journal}, 2001(2):111--125.

\bibitem[Beirlant et~al., 2010]{Beirlant2010}
Beirlant, J., Guillou, A., and Toulemonde, G. (2010).
\newblock Peaks-over-threshold modeling under random censoring.
\newblock {\em Communications in Statistics - Theory and Methods},
  39(7):1158--1179.

\bibitem[Beran, 1981]{Beran1981}
Beran, R. (1981).
\newblock Nonparametric regression with randomly censored survival data.
\newblock {\em University of California, Berkeley}.

\bibitem[Chown et~al., 2018]{Chown2018}
Chown, J., Heuchenne, C., and Van~Keilegom, I. (2018).
\newblock The nonparametric location-scale mixture cure model.
\newblock {\em TEST (under revision)}.

\bibitem[de~Haan and Ferreira, 2006]{deHaan2006}
de~Haan, L. and Ferreira, A. (2006).
\newblock {\em Extreme {V}alue {T}heory. An {I}ntroduction}.
\newblock Springer Series in Operations Research and Financial Engineering.
  Springer, New York.

\bibitem[Einmahl et~al., 2008]{Einmahl2008}
Einmahl, J.~H., Fils-Villetard, A., and Guillou, A. (2008).
\newblock Statistics of extremes under random censoring.
\newblock {\em Bernoulli}, 14(1):207--227.

\bibitem[Escobar-Bach and Van~Keilegom, 2019]{Escobar2019}
Escobar-Bach, M. and Van~Keilegom, I. (2019).
\newblock Non-parametric cure rate estimation under insufficient follow-up
  using extremes.
\newblock {\em Journal of the Royal Statistical Society - Series B (to
  appear)}.

\bibitem[Fernique, 1964]{Fernique1964}
Fernique, X. (1964).
\newblock Continuit\'{e} des processus {G}aussiens.
\newblock {\em C. R. Acad. Sci. Paris}, 258:6058--6060.

\bibitem[Fleming and Harrington, 1991]{Fleming1991}
Fleming, T.~R. and Harrington, D.~P. (1991).
\newblock {\em Counting {P}rocesses and {S}urvival {A}nalysis}.
\newblock Wiley Series in Probability and Mathematical Statistics: Applied
  Probability and Statistics. John Wiley \& Sons, Inc., New York.

\bibitem[Gin\'e and Guillou, 2002]{Gine2002}
Gin\'e, E. and Guillou, A. (2002).
\newblock Rates of strong uniform consistency for multivariate kernel density
  estimators.
\newblock {\em Annales de l'I.H.P. Probabilit\'es et Statistiques},
  38(6):907--921.

\bibitem[Gomes and Neves, 2011]{Gomes2011}
Gomes, M.~I. and Neves, M.~M. (2011).
\newblock Estimation of the extreme value index for randomly censored data.
\newblock {\em Biometrical Letters}, 48(1):1--22.

\bibitem[Gonzalez-Manteiga and Cadarso-Suarez, 1994]{Gonzalez1994}
Gonzalez-Manteiga, W. and Cadarso-Suarez, C. (1994).
\newblock Asymptotic properties of a generalized {K}aplan-{M}eier estimator
  with some applications.
\newblock {\em Journal of Nonparametric Statistics}, 4(1):65--78.

\bibitem[Kaplan and Meier, 1958]{Kaplan1958}
Kaplan, E. and Meier, P. (1958).
\newblock Nonparametric estimation from incomplete observations.
\newblock {\em Journal of the American Statistical Association}, 53:457--481.

\bibitem[L\'opez-Cheda et~al., 2017a]{Lopez2017}
L\'opez-Cheda, A., Cao, R., J\'acome, M., and Van~Keilegom, I. (2017a).
\newblock Nonparametric incidence estimation and bootstrap bandwidth selection
  in mixture cure models.
\newblock {\em Computational Statistics and Data Analysis}, 105:144--165.

\bibitem[L\'opez-Cheda et~al., 2017b]{Lopez2017b}
L\'opez-Cheda, A., J\'acome, M., and Cao, R. (2017b).
\newblock Nonparametric latency estimation for mixture cure models.
\newblock {\em TEST}, 26:353--376.

\bibitem[Maller and Zhou, 1992]{Maller1992}
Maller, R.~A. and Zhou, S. (1992).
\newblock Estimating the proportion of immunes in a censored sample.
\newblock {\em Biometrika}, 79(4):731--739.

\bibitem[Nolan and Pollard, 1987]{Nolan1987}
Nolan, D. and Pollard, D. (1987).
\newblock {$U$}-processes: rates of convergence.
\newblock {\em Annals of Statistics}, 15(2):780--799.

\bibitem[Peng and Taylor, 2014]{Peng2014}
Peng, Y. and Taylor, J. M.~G. (2014).
\newblock Cure models.
\newblock {\em Handbook of Survival Analysis, Handbooks of Modern Statistical
  Methods series, ed. J Klein, H van Houwelingen, JG Ibrahim, TH Scheike,
  chapter 6}, pages 113--134.

\bibitem[Stupfler, 2016]{Stupfler2016}
Stupfler, G. (2016).
\newblock Estimating the conditional extreme-value index under random
  right-censoring.
\newblock {\em Journal of Multivariate Analysis}, 144:1--24.

\bibitem[Van~der Vaart, 1998]{Vaart1998}
Van~der Vaart, A.~W. (1998).
\newblock {\em Asymptotic {S}tatistics}, volume~3 of {\em Cambridge Series in
  Statistical and Probabilistic Mathematics}.
\newblock Cambridge University Press, Cambridge.

\bibitem[Van~der Vaart and Wellner, 1996]{Vaart1996}
Van~der Vaart, A.~W. and Wellner, J.~A. (1996).
\newblock {\em Weak {C}onvergence and {E}mpirical {P}rocesses}.
\newblock Springer Series in Statistics. Springer-Verlag, New York.

\bibitem[Van~Keilegom and Veraverbeke, 1997]{VanKeilegom1997}
Van~Keilegom, I. and Veraverbeke, N. (1997).
\newblock Estimation and bootstrap with censored data in fixed design
  nonparametric regression.
\newblock {\em Annals of the Institute of Statistical Mathematics},
  49(3):467--491.

\bibitem[Worms and Worms, 2014]{Worms2014}
Worms, J. and Worms, R. (2014).
\newblock New estimators of the extreme value index under random right
  censoring, for heavy-tailed distributions.
\newblock {\em Extremes}, 17(2):337--358.

\bibitem[Xu and Peng, 2014]{Xu2014}
Xu, J. and Peng, Y. (2014).
\newblock Nonparametric cure rate estimation with covariates.
\newblock {\em Canadian Journal of Statistics}, 42(1):1--17.

\end{thebibliography}

\end{document}